\documentclass[12pt,leqno]{article}
\usepackage[doc,paper]{optional}
\usepackage{float}
\restylefloat{table*}
\usepackage{booktabs} 
\usepackage{url}
\usepackage[doc]{optional}
\usepackage{color}
\usepackage{float}
\usepackage{soul}
\usepackage{graphicx}
\definecolor{labelkey}{rgb}{0,0.08,0.45}
\definecolor{refkey}{rgb}{0,0.6,0.0}
\definecolor{lime}{rgb}{0.00,0.8,0.0}
\definecolor{lblue}{rgb}{0.5,0.5,0.99}
\usepackage[usenames,dvipsnames]{xcolor}
\usepackage{pgf}
\usepackage{tikz}
\usepackage{pgfplots}
\usetikzlibrary{shapes}
\usetikzlibrary{arrows}
\usetikzlibrary{trees}
\usetikzlibrary{decorations.pathreplacing}
\usetikzlibrary{calc}
\usetikzlibrary{patterns}
\usepackage{subcaption}
\usepackage[boxed]{algorithm2e}
\definecolor{light-gray}{gray}{0.65}
\usepackage{mathpazo}
\usepackage{amsmath}
\usepackage{amssymb}
\usepackage{theorem}
\usepackage{exscale,relsize}
\usepackage{pifont}
\newcommand{\ra}[1]{\renewcommand{\arraystretch}{#1}}

\newcommand{\timess}{\,{\textstyle\mathsmaller{\text{\ding{75}}}}}
\newcommand{\GX}{\ensuremath{\Gamma}}
\newcommand{\aw}{\ensuremath{\stackrel{\mathrm{e}}{\to}}}
\newcommand{\pw}{\ensuremath{\stackrel{\mathrm{p}}{\to}}}
\newcommand{\weakly}{\ensuremath{\:{\rightharpoonup}\:}}
\newcommand{\bDK}{$\bD{}\,$-Klee}
\newcommand{\fDK}{$\fD{}\,$-Klee}
\newcommand{\nnn}{\ensuremath{{n\in{\mathbb N}}}}
\newcommand{\knn}{\ensuremath{{k\in{\mathbb N}}}}
\newcommand{\thalb}{\ensuremath{\tfrac{1}{2}}}
\newcommand{\menge}[2]{\big\{{#1}~\big |~{#2}\big\}}
\newcommand{\mmenge}[2]{\bigg\{{#1}~\bigg |~{#2}\bigg\}}
\newcommand{\Menge}[2]{\left\{{#1}~\Big|~{#2}\right\}}
\newcommand{\todo}{\hookrightarrow\textsf{TO DO:}}
\newcommand{\To}{\ensuremath{\rightrightarrows}}
\newcommand{\lev}[1]{\ensuremath{\mathrm{lev}_{\leq #1}\:}}
\newcommand{\fenv}[1]%
{\ensuremath{\,\overrightarrow{\operatorname{env}}_{#1}}}
\newcommand{\benv}[1]%
{\ensuremath{\,\overleftarrow{\operatorname{env}}_{#1}}}
\newcommand{\emp}{\ensuremath{\varnothing}}
\newcommand{\infconv}{\ensuremath{\mbox{\small$\,\square\,$}}}
\newcommand{\pair}[2]{\left\langle{{#1},{#2}}\right\rangle}
\newcommand{\scal}[2]{\left\langle{#1},{#2}  \right\rangle}
\newcommand{\Tt}{\ensuremath{\mathfrak{T}}}
\newcommand{\YY}{\ensuremath{\mathcal Y}}
\newcommand{\mA}{\ensuremath{\mathcal A}}
\newcommand{\HH}{\ensuremath{\mathcal H}}
\newcommand{\exi}{\ensuremath{\exists\,}}
\newcommand{\zeroun}{\ensuremath{\left]0,1\right[}}
\newcommand{\RR}{\ensuremath{\mathbb R}}
\newcommand{\ZZ}{\ensuremath{\mathbb Z}}
\newcommand{\RP}{\ensuremath{\mathbb{R}_+}}
\newcommand{\RPX}{\ensuremath{\left[0,+\infty\right]}}
\newcommand{\RPP}{\ensuremath{\,\left]0,+\infty\right[}}
\newcommand{\RX}{\ensuremath{\,\left]-\infty,+\infty\right]}}
\newcommand{\RRX}{\ensuremath{\,\left[-\infty,+\infty\right]}}
\newcommand{\ball}{\ensuremath{\mathbb B}}
\newcommand{\KK}{\ensuremath{\mathbb K}}
\newcommand{\NN}{\ensuremath{\mathbb N}}
\newcommand{\CC}{\ensuremath{\mathbb C}}
\newcommand{\SSS}{\ensuremath{\mathbb S}}
\newcommand{\oldIDD}{\ensuremath{\operatorname{int}\operatorname{dom}f}}
\newcommand{\IDD}{\ensuremath{U}}
\newcommand{\BIDD}{\ensuremath{\mathbf U}}
\newcommand{\CDD}{\ensuremath{\overline{\operatorname{dom}}\,f}}
\newcommand{\dom}{\ensuremath{\operatorname{dom}}}
\newcommand{\argmin}{\ensuremath{\operatorname{argmin}}}
\newcommand{\argmax}{\ensuremath{\operatorname*{argmax}}}
\newcommand{\cont}{\ensuremath{\operatorname{cont}}}
\newcommand{\gr}{\ensuremath{\operatorname{gr}}}
\newcommand{\prox}{\ensuremath{\operatorname{Prox}}}
\newcommand{\intdom}{\ensuremath{\operatorname{int}\operatorname{dom}}\,}
\newcommand{\reli}{\ensuremath{\operatorname{ri}}}
\newcommand{\inte}{\ensuremath{\operatorname{int}}}
\newcommand{\deriv}{\ensuremath{\operatorname{\; d}}}
\newcommand{\rockderiv}{\ensuremath{\operatorname{\; \hat{d}}}}
\newcommand{\closu}{\ensuremath{\operatorname{cl}}}
\newcommand{\cart}{\ensuremath{\mbox{\LARGE{$\times$}}}}
\newcommand{\SC}{\ensuremath{{\mathfrak S}}}
\newcommand{\card}{\ensuremath{\operatorname{card}}}
\newcommand{\bd}{\ensuremath{\mathbf{d}}}
\newcommand{\ran}{\ensuremath{\operatorname{ran}}}
\newcommand{\conv}{\ensuremath{\operatorname{conv}}}
\newcommand{\aff}{\ensuremath{\operatorname{aff}}}
\newcommand{\spa}{\ensuremath{\operatorname{span}}}
\newcommand{\cconv}{\ensuremath{\overline{\operatorname{conv}}\,}}
\newcommand{\cdom}{\ensuremath{\overline{\operatorname{dom}}}}
\newcommand{\Fix}{\ensuremath{\operatorname{Fix}}}
\newcommand{\Id}{\ensuremath{\operatorname{Id}}}
\newcommand{\fprox}[1]{\overrightarrow{\operatorname{prox}}_{\thinspace#1}}
\newcommand{\bprox}[1]{\overleftarrow{\operatorname{prox}}_{\thinspace#1}}
\newcommand{\Hess}{\ensuremath{\nabla^2\!}}
\newcommand{\fproj}[1]{\overrightarrow{P\thinspace}_%
{\negthinspace\negthinspace #1}}
\newcommand{\ffproj}[1]{\overrightarrow{Q\thinspace}_%
{\negthinspace\negthinspace #1}}
\newcommand{\bproj}[1]{\overleftarrow{\thinspace P\thinspace}_%
{\negthinspace\negthinspace #1}}
\newcommand{\bfproj}[1]{\overleftarrow{\thinspace Q\thinspace}_%
{\negthinspace\negthinspace #1}}
\newcommand{\proj}[1]{{\thinspace P\thinspace}_%
{\negthinspace\negthinspace #1}}
\newcommand{\fD}[1]{\overrightarrow{D\thinspace}_%
{\negthinspace\negthinspace #1}}
\newcommand{\ffD}[1]{\overrightarrow{F\thinspace}_%
{\negthinspace\negthinspace #1}}
\newcommand{\ffDbz}{\protect\overrightarrow{F\thinspace}_%
{\negthinspace\negthinspace C}(\bz)}
\newcommand{\bfD}[1]{\overleftarrow{\thinspace F\thinspace}_%
{\negthinspace\negthinspace #1}}

\newcommand{\jj}{\ensuremath{\,\mathfrak{q}}}

\newcommand{\anfang}{\ensuremath{{{v}}}}
\newcommand{\ziel}{\ensuremath{{{w}}}}
\newcommand{\minf}{\ensuremath{-\infty}}
\newcommand{\pinf}{\ensuremath{+\infty}}
\newcommand{\bx}{\ensuremath{\mathbf{x}}}
\newcommand{\ba}{\ensuremath{\mathbf{a}}}
\newcommand{\bX}{\ensuremath{{\mathbf{X}}}}
\newcommand{\bR}{\ensuremath{{\mathbf{R}}}}
\newcommand{\bT}{\ensuremath{{\mathbf{T}}}}
\newcommand{\wT}{\ensuremath{{\widetilde{\mathbf{T}}}}}
\newcommand{\bN}{\ensuremath{{\mathbf{N}}}}
\newcommand{\bL}{\ensuremath{{\mathbf{L}}}}
\newcommand{\bJ}{\ensuremath{{\mathbf{J}}}}
\newcommand{\bS}{\ensuremath{{\mathbf{S}}}}
\newcommand{\bD}{\ensuremath{{\mathbf{D}}}}
\newcommand{\bA}{\ensuremath{{{\mathbf{A}}}}}
\newcommand{\bQ}{\ensuremath{{{\mathbf{Q}}}}}
\newcommand{\bId}{\ensuremath{{\mathbf{Id}}}}
\newcommand{\bB}{\ensuremath{{\mathbf{B}}}}
\newcommand{\bc}{\ensuremath{\mathbf{c}}}
\newcommand{\bee}{\ensuremath{\mathbf{e}}}
\newcommand{\bC}{\ensuremath{\mathbf{C}}}
\newcommand{\by}{\ensuremath{\mathbf{y}}}
\newcommand{\bz}{\ensuremath{\mathbf{z}}}
\newcommand{\bg}{\ensuremath{\mathbf{g}}}
\newcommand{\sm}{\ensuremath{\mathbb{S}^{N\times N}}}
\newcommand{\smp}{\ensuremath{\mathbb{S}^{N\times N}_{+}}}
\newcommand{\smpp}{\ensuremath{\mathbb{S}^{N\times N}_{++}}}
\newcommand{\bh}{\ensuremath{\mathbf{h}}}
\newcommand{\Tbar}{\ensuremath{\overline{T}}}
\newcommand{\pfl}{\ensuremath{p({\fettf,\fettla})}}
\newcommand{\pflm}{\ensuremath{p_\mu({\fettf,\fettla)}}}
\newcommand{\fettf}{\ensuremath{\boldsymbol{f}}}
\newcommand{\fettla}{\ensuremath{\boldsymbol{\lambda}}}
\newcommand{\ve}{\ensuremath{{{\varepsilon}}}}

\DeclareMathOperator{\sgn}{sgn}

{\begin{list}{}{\renewcommand{\makelabel}[1]{\textrm{##1~}\hfil}%
\settowidth{\labelwidth}{\textrm{#1~}}%
\setlength{\leftmargin}{\labelwidth+\labelsep}}}
{\end{list}}
\newtheorem{theorem}{Theorem}[section]

\newtheorem{proposition}[theorem]{Proposition}

\theoremstyle{plain}{\theorembodyfont{\rmfamily}
}
\theoremstyle{plain}{\theorembodyfont{\rmfamily}
}
\theoremstyle{plain}{\theorembodyfont{\rmfamily}
\newtheorem{algo}[theorem]{Algorithm}}
\theoremstyle{plain}{\theorembodyfont{\rmfamily}
\newtheorem{example}[theorem]{Example}}

\theoremstyle{plain}{\theorembodyfont{\rmfamily}
}
\theoremstyle{plain}{\theorembodyfont{\rmfamily}
\newtheorem{remark}[theorem]{Remark}}
\def\proof{\noindent{\it Proof}. \ignorespaces}
\def\endproof{\ensuremath{\hfill \quad \blacksquare}}
\renewcommand\theenumi{(\roman{enumi})}
\renewcommand\theenumii{(\alph{enumii})}
\renewcommand{\labelenumi}{\rm (\roman{enumi})}
\renewcommand{\labelenumii}{\rm (\alph{enumii})}

\newcommand{\boxedeqn}[1]{%
    \[\fbox{%
        \addtolength{\linewidth}{-2\fboxsep}%
        \addtolength{\linewidth}{-2\fboxrule}%
        \begin{minipage}{\linewidth}%
        \begin{equation}#1\\[+4mm]\end{equation}%
        \end{minipage}%
      }\]%
  }

\newcommand{\algJA}{\text{\texttt{alg}}(J_A)}
\newcommand{\algJR}{{\text{\texttt{alg}}(\bJ\circ\bR)}}
\newcommand{\algT}{\text{\texttt{alg}}(\bT)}

\allowdisplaybreaks

\begin{document}

\title{\textrm{Projection Methods: Swiss Army Knives for Solving
Feasibility and Best Approximation Problems with Halfspaces}}

\author{
Heinz H.\ Bauschke\thanks{
Mathematics, University
of British Columbia,
Kelowna, B.C.\ V1V~1V7, Canada.
E-mail: \texttt{heinz.bauschke@ubc.ca}.}
~and~
Valentin R.\ Koch\thanks{
Information Modeling \& Platform Products Group (IPG),
Autodesk, Inc.
E-mail:  \texttt{valentin.koch@autodesk.com}.}
}

\date{January 18, 2013}

\maketitle

\vskip 8mm

\begin{abstract} \noindent
We model a problem motivated by road design as a
feasibility problem. Projections onto
the constraint sets are obtained, and projection methods
for solving the feasibility problem are studied.
We present results of numerical experiments which 
demonstrate the efficacy of projection methods even
for challenging nonconvex problems. 
\end{abstract}

{\small


}

\section{Introduction and motivation}

\subsection{The abstract formulation of the problem}
\label{ss:abstract}

Throughout this paper, we assume that $n\in\{2,3,\ldots\}$ and that 
\boxedeqn{
\text{$X=\RR^n$ 
with standard inner product $\scal{\cdot}{\cdot}$
and induced norm $\|\cdot\|$.
}
}
We also assume we are given $n$ strictly increasing breakpoints on the real
line: 
\boxedeqn{\label{e:xticks}
t = (t_1,\ldots,t_n)\in X
\quad\text{such that $t_1 < \cdots < t_n$.}
}
Our goal is to 
\begin{subequations}
\label{e:ubergoal}
\begin{equation}
\label{e:abstractgoal}
\text{find a vector $x=(x_1,\ldots,x_n)\in X$}
\end{equation}
such that 
\begin{equation}
\label{e:abstractconstraints}
\text{ $(t_1,x_1),\ldots,(t_n,x_n)$ satisfies a given set of constraints.}
\end{equation}
\end{subequations}
Note that for every $x\in X$, 
the pair $(t,x)\in X\times X$ 
induces a corresponding piecewise linear function
or \emph{linear spline} (see \cite{deBoor} and \cite{Schumaker}) 
\begin{equation}
f_{(t,x)}\colon [t_1,t_n] \to \RR\colon
\tau \mapsto x_{i} + (x_{i+1}-x_i)\frac{\tau -t_i}{t_{i+1}-t_i},
\end{equation}
where $\tau \in [t_i,t_{i+1}]$ and $i\in\{1,\ldots,n-1\}$,
which passes through the points $\menge{(t_i,x_i)\in\RR^2}{i\in I}$.

The set of constraints mentioned in \eqref{e:abstractconstraints} will
involve the function $f_{(t,x)}$. 
Let us list several types of constraints which are
motivated in Section~\ref{ss:Applications} below:

\begin{itemize}
\item 
\textbf{interpolation constraints:}
For a given subset $I$ of
$\{1,\ldots,n\}$, 
the entries $x_i$ are prescribed:
$(\forall i\in I)$ $f_{(t,x)}(t_i) = y_i$. 
This is an \emph{interpolation problem} for
the points $\menge{(t_i,y_i)\in\RR^2}{i\in I}$. 
\item
\textbf{slope constraints:}
For a given subset $I$ of 
$\{1,\ldots,n-1\}$ and for every $i\in I$, 
the \emph{slope} 
\begin{equation}
s_i := \frac{x_{i+1}-x_i}{t_{i+1}-t_i}
\end{equation}
of $f_{(t,x)}|_{[t_i,t_{i+1}]}$ must lie in a given subset of $\RR$. 
\item
\textbf{curvature constraints:}
For a given subset $I$ of 
$\{1,\ldots,n-2\}$ and for every $i\in I$, 
$|s_{i+1}-s_i|$, 
the distance between the slopes $s_i$ and $s_{i+1}$ of two adjacent
intervals $[t_i,t_{i+1}]$ and $[t_{i+1},t_{i+2}]$, 
must lie in a given subset of $\RR$.
\end{itemize}

\subsection{A concrete instance in road design}
\label{ss:Applications}

The problem~\eqref{e:ubergoal} introduced above
and its solutions have several direct applications in
engineering and computer-assisted design. 
For instance, an engineer may want to verify the feasibility of a
design, or adapt the design according to the constraints. 
Examples drawn from Computer-Assisted Design (CAD)
include designs for roadway profiles,
pipe layouts, fuel lines in automotive designs such as cars and
airplanes, overhead power lines, chairlifts, cable cars, and duct
networks.

Our primary motivation for this work is 
automated design of \emph{road alignments}. 
A road alignment is represented 
by the centerline of the road, which is idealized as 
a (generally) nonlinear, smooth curve in $\RR^3$.
To facilitate construction drawings, civil engineers reduce the
three-dimensional road design to two two-dimensional parts, horizontal and
vertical.

\begin{figure}[H]
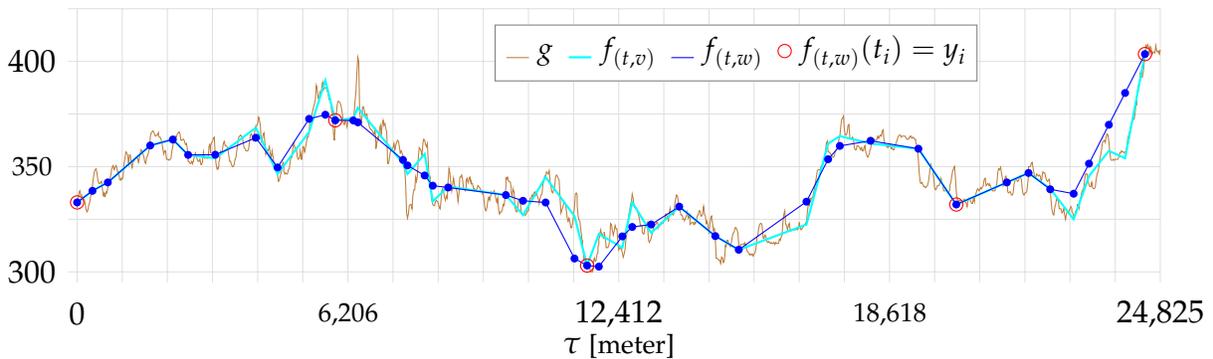
 
\centering

\caption{Vertical profiles of ground profile $g$, initial road
design $f_{(t,\anfang)}$, and final road design $f_{(t,\ziel)}$ for
a highway design near Kelowna, B.C., Canada.
The solution was found
with the ParDyk algorithm (Algorithm~\ref{a:ParDyk} below)
using a design speed of 130 km/h, and a maximum slope of 4\%. 
(These engineering constraints translate into a specific set of
slope and interpolation constraints.)
}
\label{fig:road}
\end{figure}

The horizontal alignment is the plan (or map) view of the road. 
In the vertical view, the \emph{ground profile} 
$g \colon [t_1,t_n] \to \RR$ shows
the elevation values of the existing ground along the centerline
(see the brown curve in Figure~\ref{fig:road}). 
Since earthwork operations
such as cut and fills are expensive items in road construction, a
civil engineer tries to find a road profile represented by 
a linear spline $f_{(t,x)}$) that follows $g$ as closely as possible.

Design constraints imposed on $f_{(t,x)}$ by 
the engineer or by 
civil design standards such as those published by the 
\emph{American Association of State Highway and Transportation
Officials (AASHTO)} \cite{AASHTO} 
include the following:
\begin{itemize}
\item 
At a certain station $t_i$, the engineer may have to fix the
elevation $x_i$ to allow for construction of an intersection with
an existing road that crosses the new road at $t_i$. This 
corresponds to the (mathematically \emph{affine}) 
interpolation constraint mentioned in Section~\ref{ss:abstract}. 
\item For safety reasons and to ensure good
traffic flow, AASHTO requires that the slope between two stations
$t_i$ and $t_{i+1}$ is bounded above and below. These are the
(mathematically \emph{convex}) slope constraints of
Section~\ref{ss:abstract}. 
\item In a road profile, engineers
fit a (usually parabolic) curve at the point of intersection
of two line segments. The curvature depends on the grade change of
the line segments and influences the vertical acceleration of a
vehicle, as well as the stopping sight distance. AASHTO requires
bounds on the slope change. This corresponds to the 
(mathematically \emph{convex}) curvature
constraint in Section~\ref{ss:abstract}. 
\item In some cases, the engineer requires a minimum
drainage grade to allow flow and avoid 
catchment of storm water. 
These (mathematically challenging \emph{nonconvex})
slope constraints are discussed in Section~\ref{sec:slope-noncon}
below.
\end{itemize}
We denote the starting spline for a road profile 
(see the cyan curve in Figure \ref{fig:road}) by 
\begin{equation}
f_{(t,\anfang)}, \quad
\text{where $t = (t_1,\ldots,t_n)\in X$ and
$\anfang = (\anfang_1,\ldots,\anfang_n)\in X$.}
\end{equation}
In practice, the starting spline could simply be the 
connected line segments for the
interpolation constraint, or it could be generated from the ground
profile $g$ by using Bellman's method 
(see \cite{Bellman} and \cite{Koch} for details).
In either case, we assume that $t$ is given and fixed, 
and that we need to decide whether or not $f_{(t,\anfang)}$
is feasible with respect to  the aforementioned constraints. 
If $\anfang$ leads to a feasible spline, then we are done;
otherwise, we wish 
to find $\ziel\in X$ such that the new road spline $f_{(t,\ziel)}$ 
(see the blue curve in Figure~\ref{fig:road}) 
satisfies the design constraints.
Ideally, $f_{(t,\ziel)}$ is close to the ground profile
represented by ${(t,\anfang)}$. 
Finally, if there is no $\ziel\in X$ making the problem feasible, 
then we would like to detect this through some suitable measure. 

\subsection{Main results and organization of the paper}

We now summarize the main contributions of this paper.

\begin{itemize}
\item In principle, there are numerous constraints to deal with for the
problem \eqref{e:ubergoal} in the context of road design.
Fortunately, the constraints have a lot of structure we can take
advantage of and we demonstrate that 
\emph{the constraints parallelize which allows to reduce
the problem to six constraint sets} 
(see Section~\ref{ss:constraintssummary})
each of which admits a closed form projection
formula (see \eqref{e:PY}, \eqref{e:PS}, and \eqref{e:Pcurv}, 
\eqref{e:simpleslope+}). 
\item We provide a selection of \emph{state-of-the-art projection
methods, 
superiorization algorithms, and best approximation algorithms} 
(see Sections~\ref{sec:FPM} and \ref{sec:ba}), 
and adapt them to the road design problem. 
\item
We \emph{present various observations on the algorithms and their 
relationships} 
(see Remark~\ref{r:DRvsADMM}, Remark~\ref{r:DykvsADMM}, 
Remark~\ref{r:doublyaffine}, Remark~\ref{r:notequal}
and Example~\ref{ex:cycling}.)
\item 
We report on broad \emph{numerical experiments} (see Section~\ref{sec:numeresul}) 
introducing for the first time \emph{performance profiles} for projection
methods. 
\item Based on the numerical experiments, we recommend CycP+, which is
an intrepid form of the method of cyclic projections, as an overall good
choice for solving feasibility and best approximation problems. 
\end{itemize}

The remainder of the paper is organized as follows.
Section~\ref{sec:CPO} contains a detailed analysis of the projection
operators encountered in the road design problem. We take advantage
of aggregating constraints and derive simple formulas. 
Projection methods for feasibility problems are reviewed in
Section~\ref{sec:FPM}. Because we are working with more than two constraint
sets, we adapt to this situation by working in a product space if
needed. If more than just a feasible road design is desired, then
the engineer has to consider optimization algorithms. 
We review two types of such methods
(superiorization and best approximation algorithms) and adapt them to our problem
in Section~\ref{ss:superior} and Section~\ref{sec:ba} respectively.
Nonconvex constraints are investigated in
Section~\ref{sec:nonconvexity}. 
We report on numerical experiments in Section~\ref{sec:numeresul}.
The final Section~\ref{sec:conclusion} concludes the paper. 

For notation and general references on the mathematics underlying
projection methods, we refer the reader to the books
\cite{BC2011}, \cite{Cegielski}, \cite{CZ}, \cite{Deutsch},
\cite{GK}, and \cite{GR}.

\section{Constraints and projection operators}

\label{sec:CPO}

\subsection{The projection onto a general convex set}

In this section, we make the constraints encountered
in road design mathematically precise.
Almost all of these constraints turn out to lead to sets that 
are \emph{convex} and closed.
Recall that a set $C$ is \emph{convex} if it contains all line segments 
between each pair taken from $C$: 
\begin{equation}
(\forall c_0\in C)(\forall c_1\in C)(\forall \lambda\in[0,1])
\quad (1-\lambda)c_0 + \lambda c_1\in C.
\end{equation}
If $C$ is a nonempty closed convex subset of $X$, 
then for every $x\in X$, the optimization problem
\begin{equation}
\label{e:distance}
d(x,C) := \min_{c\in C} \|x-c\|,
\end{equation}
which concerns the computation of the \emph{distance} $d(x,C)$ from
$x$ to the set $C$, 
has a \emph{unique} solution, denoted by $P_Cx$
and called the \emph{projection} of $x$ onto $C$. 
The vector $P_Cx$ is characterized by two properties, namely 
\begin{equation}
\label{e:kolmogorov}
P_Cx\in C
\quad\text{and}\quad
(\forall c\in C)\quad \scal{c-P_Cx}{x-P_Cx} \leq 0.
\end{equation}
(For a proof, see, e.g., \cite[Theorem~3.14]{BC2011}). 
The induced operator $P_C\colon X\to C$ 
is the \emph{projection operator} or \emph{projector} onto $C$. 
There are several examples that allow us to write down the projector in closed
form or to approximate $P_C$ provided $C$ is the intersection of finitely
many simple closed convex sets (see, e.g., \cite[Chapters~28 and 29]{BC2011}.)
In the road design application, it is fortunately 
possible to obtain explicit formulas; in the following subsections, 
we will make these formulas as convenient as possible
for software implementation. 
Let us do this for each of the three types of constraints.

\subsection{Interpolation constraints}

We assume that $I$ is a set such that 
\begin{equation}
\{1,n\} \subseteq I \subseteq \{1,2,\ldots,n\},
\quad\text{and}\quad
{y} = ({y}_i)\in\RR^I
\end{equation}
is given.
Set 
\begin{equation}
\label{e:Y}
Y := \menge{x=(x_1,\ldots,x_n)\in X}{(\forall i\in I)\;\; x_i = {y}_i}.
\end{equation}
The closed set $Y$ is an \emph{affine subspace}, i.e.,
$(\forall y\in Y)(\forall z\in Y)(\forall \lambda\in\RR)$
$(1-\lambda)y + \lambda z\in Y$; in particular, $Y$ is convex. 
For convenience, 
we record the explicit formula for the projection onto $Y$.

\begin{proposition}[interpolation constraint projector]
\label{p:PY}
The projector onto $Y$ is given by 
\begin{subequations}
\label{e:PY}
\begin{align}
P_Y\colon X&\to X\\
(x_1,x_2,\ldots,x_n)&\mapsto (c_1,c_2,\ldots,c_n),\quad
\text{where}\quad
c_i = \begin{cases}
{y}_i, &\text{if $i\in I$;}\\
x_i, &\text{if $i\in \{1,2,\ldots,n\}\smallsetminus I$.}
\end{cases}
\end{align}
\end{subequations}
\end{proposition}

\subsection{Slope constraints}

\subsubsection{A useful special case}

\label{sss:slopespecial}

We start with a simple special case that will also be useful
in handling our general slope constraints.
To this end, 
let $i\in \{1,2,\ldots,n-1\}$. The constraint set $S_i$ imposes 
that 
the absolute value of the slope $f_{(t,x)}$ for the interval $[t_i,t_{i+1}]$
is bounded above, i.e., 
there exists $\alpha_i\geq 0$ such that 
\begin{equation}
\label{e:Si}
S_i := \menge{x=(x_1,\ldots,x_n)\in\RR^n}{
|x_{i+1}-x_i| \leq \alpha_i}.
\end{equation}
Indeed, 
if the actual maximum absolute slope is $\sigma_i \geq 0$,
then, setting $\alpha_i = \sigma_i|t_{i+1}-t_i|$, we see that
\begin{subequations}
\begin{align}
\left| \frac{x_{i+1}-x_i}{t_{i+1}-t_i}\right| \leq \sigma_i
&\Leftrightarrow
|x_{i+1}-x_i| \leq \alpha_i
\Leftrightarrow
\pm(x_{i+1}-x_i) \leq \alpha_i. \\
&\Leftrightarrow
\pm\scal{e_{i+1}-e_i}{x} \leq \alpha_i\\
&\Leftrightarrow
-\alpha_i\leq \scal{e_{i+1}-e_i}{x} \leq \alpha_i,
\end{align}
\end{subequations}
where $e_i$ and $e_{i+1}$ denote the standard unit vectors in $X$ (with the
number $1$ in position $i$ and $i+1$, respectively, and zeros
elsewhere). 
The last characterization reveals that $S_i$ is the intersection
of two halfspaces whose boundary hyperplanes are orthogonal to 
$e_{i+1}-e_i$. 
In particular, $S_i$ is a closed convex subset of $X$. 
Using e.g., \cite[Example~28.17]{BC2011}, we obtain for
every $x\in X$, 
\begin{equation}
P_{S_i}x = 
\begin{cases}
x + \displaystyle \frac{-\alpha_i-\scal{e_{i+1}-e_i}{x}}{\|e_{i+1}-e_i\|^2}(e_{i+1}-e_i),
&\text{if $\scal{e_{i+1}-e_i}{x} < -\alpha_i$;}\\
x, &\text{if $-\alpha_i\leq \scal{e_{i+1}-e_i}{x} \leq \alpha_i$;}\\
x + \displaystyle \frac{\alpha_i-\scal{e_{i+1}-e_i}{x}}{\|e_{i+1}-e_i\|^2}(e_{i+1}-e_i),
&\text{if $\scal{e_{i+1}-e_i}{x} > \alpha_i$.}\\
\end{cases}
\end{equation}
This formula shows that $(P_{S_i}x)_j = x_j$ for every $j\in
I\smallsetminus\{i,i+1\}$. 
Thus, the only entries that possibly
change after executing $P_{S_i}$ are in positions $i$ and $i+1$;
after some simplification, we obtain for these entries the formula
\begin{equation}
\label{e:simpleslope}
\big( (P_{S_i}x)_i, (P_{S_i}x)_{i+1}\big) = 
\begin{cases}
\thalb(x_i+x_{i+1}+\alpha_i,x_i+x_{i+1}-\alpha_i), &\text{if
$x_{i}-x_{i+1}>\alpha_i$;}\\
(x_i,x_{i+1}), &\text{if $|x_{i+1}-x_i|\leq\alpha_i$;}\\
\thalb(x_i+x_{i+1}-\alpha_i,x_i+x_{i+1}+\alpha_i), &\text{if
$x_{i+1}-x_i>\alpha_i$.}
\end{cases}
\end{equation}
We note that 
\begin{equation}
\label{e:slopetrouble}
\text{$x\notin S_i\;\;\Rightarrow\;\;
(P_{S_i}x)_i \neq x_i$ and $(P_{S_i}x)_{i+1}\neq x_{i+1}$;}
\end{equation}
furthermore, if $\alpha_i = \pinf$, i.e.,
no slope constraint, then 
\eqref{e:simpleslope} is valid as well. 

\subsubsection{The general case}

\label{sss:slopegeneral}

Now we turn to the general case. 
We assume the existence of a vector $a = (\alpha_i)\in\RR_+^{n-1}$
such that the constraint set is
\begin{equation}
\label{e:Siall}
S = 
\bigcap_{i\in\{1,\ldots,n-1\}}S_i = 
\menge{(x_1,\ldots,x_n)\in\RR^n}{
(\forall i\in\{1,\ldots,n-1\})\;\;|x_{i+1}-x_i| \leq \alpha_i},
\end{equation}
where $S_i$ is defined in \eqref{e:Si}. 
While we obtained an explicit formula to deal with 
a single slope constraint (see \eqref{e:simpleslope}) 
we are unaware of a corresponding formula for $P_S$.
Furthermore, since $P_{S_i}$ possibly modifies the vector in
positions $i$ and $i+1$ (but not elsewhere), 
we \emph{cannot} use \eqref{e:simpleslope} for the sets $S_i$ and
$S_{i+1}$ concurrently because their projections
possibly modify positions $(i,i+1)$ and $(i+1,i+2)$ (see
\eqref{e:slopetrouble}), but not necessarily in a consistent manner
at position $i+1$!
However, 
by combining the $n-1$ slope constraints according to parity of indices,
i.e., by setting 
\begin{equation}
\label{e:S}
S_{\text{\rm even}} := 
\bigcap_{i\in\{1,\ldots,n-1\}\cap(2\NN)}S_i,
\quad\text{and}\quad
S_{\text{\rm odd}} := 
\bigcap_{i\in\{1,\ldots,n-1\}\cap(1+2\NN)}S_i,
\end{equation}
we see that
\begin{equation}
S = S_{\text{\rm even}} \cap S_{\text{\rm odd}} 
\end{equation}
can be written of the intersection of just \emph{two} constraint sets!
Furthermore, \eqref{e:simpleslope} yields the 
\emph{fully parallel} update formulas:

\begin{proposition}[convex slope constraint projector]
\label{p:PS}
For every $x\in X$, the projectors onto $S_{\text{even}}$
and $S_{\text{odd}}$ are given by 
\begin{subequations}
\label{e:PS}
\begin{equation}
\label{e:PSeven}
P_{S_{\text{\rm even}}}x = 
\big(x_1,(P_{S_2}x)_2,(P_{S_2}x)_3,(P_{S_4}x)_4,(P_{S_4}x)_5,\ldots\big)\in
X,
\end{equation}
where the last entry in \eqref{e:PSeven} is $x_n$ if $n$ is even; 
and 
\begin{equation}
\label{e:PSodd}
P_{S_{\text{\rm odd}}}x = 
\big((P_{S_1}x)_1,(P_{S_1}x)_2,(P_{S_3}x)_3,(P_{S_3}x)_4,\ldots\big)\in X,
\end{equation}
\end{subequations}
where the last entry in \eqref{e:PSodd} is $x_n$ if $n$ is odd. 
\end{proposition}

The constraints making up the aggregated slope constraint are very
special polyhedra, namely \emph{``stripes''}, i.e., 
the intersection of two halfspaces with opposing
normal vectors. 
This motivates the technique, 
which originates in \cite{HC} (see also \cite{CCCDH}, \cite{CCP}, and
\cite{Herman}), 
of not just projecting onto these sets but rather
\emph{inside} them: either we reflect \emph{into} the strip or
(if we are too distant from the strip) we 
jump to the corresponding midpoint of the strip. 
Let us record the formula for this operator.

\begin{proposition}[intrepid slope constraint projectors]
\label{p:iPS}
The intrepid counterpart of
\eqref{e:simpleslope}, is\footnote{\ Here $\sgn$ denotes the signum
function.} 
\begin{equation}
\label{e:aggressivesimpleslope}
(x_i,x_{i+1})\mapsto 
\begin{cases}
\thalb(x_i+x_{i+1},x_i+x_{i+1}), &\text{if
$|x_{i}-x_{i+1}|>2\alpha_i$;}\\
(x_i,x_{i+1}), &\text{if $|x_{i+1}-x_i|\leq\alpha_i$;}\\
(x_{i+1},x_i)+\sgn(x_{i+1}-x_i)(-\alpha_i,\alpha_i), &\text{if
$\alpha_i<|x_{i}-x_{i+1}|<2\alpha_i$.}
\end{cases}
\end{equation}
These operators lead to the intrepid counterpart of \eqref{e:PS}. 
\end{proposition}

\subsection{Curvature constraints}

\subsubsection{A useful special case}

Again, let us start with a simple special case that will also be useful
in handling our general curvature constraints.
To this end, 
let $i\in \{1,\ldots,n-2\}$. The constraint requires that 
the difference of consecutive slopes is bounded above and below.
Hence there exists $\gamma_i\in\RR$ and $\delta_i\in\RR$ such that 
$\gamma_i\geq\delta_i$ and 
\begin{equation}
\label{e:Ci}
C_i := \mmenge{x=(x_1,\ldots,x_n)\in\RR^n}{\gamma_i \geq
\frac{x_{i+2}-x_{i+1}}{t_{i+2}-t_{i+1}} - 
\frac{x_{i+1}-x_{i}}{t_{i+1}-t_{i}} \geq \delta_i}.
\end{equation}
Set 
\begin{equation}
\tau_i := t_{i+1}-t_i, \;\;
\tau_{i+1} := t_{i+2}-t_{i+1},\;\;\text{and}\;\;
u_i := \tau_{i+1}e_i-(\tau_{i}+\tau_{i+1})e_{i+1}+\tau_ie_{i+2}.
\end{equation}
Then $\tau_i>0$, $\tau_{i+1}>0$, and for every $x\in X$, 
\begin{equation}
x\in C_i
\quad\Leftrightarrow\quad
\gamma_i\tau_{i+1}\tau_{i}\geq
\scal{u_i}{x}\geq
\delta_i\tau_{i+1}\tau_{i}.
\end{equation}
Again, we see that $C_i$ is a closed \emph{convex} polyhedron and by 
e.g.\ \cite[Example~28.17]{BC2011}, we obtain
\begin{equation}
\label{e:usingagain}
P_{C_i}\colon x\mapsto 
\begin{cases}
x + \displaystyle
\frac{\delta_i\tau_i\tau_{i+1}-\scal{u_i}{x}}{\|u_i\|^2} u_i,
&\text{if $\scal{u_i}{x} < \delta_i\tau_i\tau_{i+1}$;}\\
x, &\text{if $\delta_i\tau_i\tau_{i+1}\leq\scal{u_i}{x}\leq 
\gamma_i\tau_i\tau_{i+1}$;}\\
x + \displaystyle
\frac{\gamma_i\tau_i\tau_{i+1}-\scal{u_i}{x}}{\|u_i\|^2}u_i,
&\text{if $\gamma_i\tau_i\tau_{i+1}<\scal{u_i}{x} $.}
\end{cases}
\end{equation}
Since $u_i\in\spa\{e_i,e_{i+1},e_{i+2}\}$, it follows that
\begin{equation}
\label{e:3lap}
\menge{j\in\{1,\ldots,n\}}{x_j \neq (P_{C_i}x)_j}
\subseteq \{i,i+1,i+2\}.
\end{equation}

\subsubsection{The general case}

Now we turn to the general curvature constraints. 
We assume the existence of a vector $c = (\gamma_i)\in\RR^{n-2}$
and $d = (\delta_i)\in\RR^{n-2}$ 
such that the constraint set is
\begin{subequations}
\begin{align}
C &= 
\bigcap_{i\in\{1,\ldots,n-2\}}C_i\\
&= 
\mmenge{(x_1,\ldots,x_n)\in X}{(\forall i\in\{1,\ldots,n-2\})\;\gamma_i \geq
\frac{x_{i+2}-x_{i+1}}{t_{i+2}-t_{i+1}} - 
\frac{x_{i+1}-x_{i}}{t_{i+1}-t_{i}} \geq \delta_i}.
\end{align}
\end{subequations}
Because of \eqref{e:3lap}, we can and do aggregate
these $n-2$ constraints into \emph{three} constraint sets that allow
projections in closed form. To this end, 
we set
\begin{equation}
\label{e:Ccurv}
(\forall j\in\{1,2,3\})\quad
C_{[j]} := \bigcap_{i\in\{1,\ldots,n-2\}\cap(j+3\NN)} C_i
\end{equation}
so that
\begin{equation}
C = C_{[1]} \cap C_{[2]} \cap C_{[3]}.
\end{equation}

Combined with \eqref{e:usingagain}, we obtain the following.

\begin{proposition}[curvature constraint projector]
\label{p:Pcurv}
For every $x\in X$, 
the projectors onto $C_{[1]}$, $C_{[2]}$, and $C_{[3]}$ are given 
\begin{subequations}
\label{e:Pcurv}
\begin{align}
P_{C_{[1]}}x &= \big((P_{C_1}x)_1,(P_{C_1}x)_2,(P_{C_1}x)_3,
(P_{C_4}x)_4,(P_{C_4}x)_5,(P_{C_4}x)_6,\ldots \big),\\
P_{C_{[2]}}x &= \big(x_1,(P_{C_2}x)_2,(P_{C_2}x)_3,(P_{C_2}x)_4,
(P_{C_5}x)_5,(P_{C_5}x)_6,(P_{C_5}x)_7,\ldots \big),\\
P_{C_{[3]}}x &= \big(x_1,x_2,(P_{C_3}x)_3,(P_{C_3}x)_4,(P_{C_3}x)_5,
(P_{C_6}x)_6,(P_{C_6}x)_7,(P_{C_6}x)_8,\ldots \big),
\end{align}
\end{subequations}
respectively, where each $P_{C_i}$ is given by \eqref{e:usingagain}. 
\end{proposition}

We now record the intrepid counterpart.

\begin{proposition}[intrepid curvature constraint projectors]
\label{p:iPcurv}
The intrepid counterpart of \eqref{e:usingagain} is 
\begin{equation}
\label{e:aggcurpro}
x\mapsto 
\begin{cases}
x - \displaystyle
\frac{\scal{u_i}{x}-\tau_i\tau_{i+1}(\delta_i+\gamma_i)/2}{\|u_i\|^2} u_i,
&\text{if $\scal{u_i}{x} <
\tau_i\tau_{i+1}(3\delta_i-\gamma_i)/2$;}\\[+5mm]
x - 2 \displaystyle
\frac{\scal{u_i}{x}-\delta_i\tau_i\tau_{i+1}}{\|u_i\|^2} u_i,
&\text{if
$\tau_i\tau_{i+1}(3\delta_i-\gamma_i)/2\leq\scal{u_i}{x}<\delta_i\tau_i\tau_{i+1}$;}\\[+5mm]
x, &\text{if $\delta_i\tau_i\tau_{i+1}\leq\scal{u_i}{x}\leq 
\gamma_i\tau_i\tau_{i+1}$;}\\[+5mm]
x - 2 \displaystyle
\frac{\scal{u_i}{x}-\gamma_i\tau_i\tau_{i+1}}{\|u_i\|^2} u_i,
&\text{if
$\gamma_i\tau_i\tau_{i+1}<\scal{u_i}{x}\leq
\tau_i\tau_{i+1}(3\gamma_i-\delta_i)/2$;}\\[+5mm]
x - \displaystyle
\frac{\scal{u_i}{x}-\tau_i\tau_{i+1}(\delta_i+\gamma_i)/2}{\|u_i\|^2} u_i,
&\text{if $\tau_i\tau_{i+1}(3\gamma_i-\delta_i)/2<\scal{u_i}{x} $.}
\end{cases}
\end{equation}
These counterparts induce the intrepid variant of
\eqref{e:Pcurv}. 
\end{proposition}

\subsection{The convex feasibility problem}

\label{ss:constraintssummary}

The \emph{convex feasibility problem} motivated by road design
is to find a point in 
\boxedeqn{
\label{e:www}
C := Y \cap S_{\text{\rm even}}\cap S_{\text{\rm odd}}
\cap  C_{[1]} \cap C_{[2]} \cap C_{[3]},
}
where the sets on the right side are defined
in \eqref{e:Y}, \eqref{e:S}, and \eqref{e:Ccurv}. 
Note that we have explicit formulas available
for all six projectors 
(see Propositions~\ref{p:PY}, \ref{p:PS}, and \ref{p:Pcurv}).
Moreover, we have intrepid variants of the projectors onto
the five\footnote{\label{fn:PY}\ The intrepid variant of $P_Y$ is $P_Y$ itself
because the interior of the affine subspace $Y$ is empty.} constraint sets
$S_{\text{\rm even}}$, $S_{\text{\rm odd}}$, 
$C_{[1]}$, $C_{[2]}$, and $C_{[3]}$ 
(see Propositions~\ref{p:iPS} and \ref{p:iPcurv}).

\section{Feasibility problems and projection methods}

\label{sec:FPM}

In this section, we present a selection of classic projection methods
for solving feasibility problems. To this end, 
let $m\in\{2,3,\ldots\}$ and
$C_1,\ldots,C_m$ be nonempty closed convex subsets of $X$.

\subsection{The $m$-set feasibility problem and its reduction to two
sets}

\label{ss:mto2}

Our aim is to 
\begin{equation}
\label{e:cfp}
\text{find $x\in C := C_1\cap\cdots\cap C_m\neq\varnothing$,}
\end{equation}
or, even more ambitiously, to compute
\begin{equation}
\label{e:130111a}
P_{C}\anfang
\end{equation}
for some given point $\anfang\in X$, or some intermediate
between feasibility and best approximation. 
Of course, we have the concrete scenario \eqref{e:www} of
Section~\ref{ss:constraintssummary}, where $m=6$, in mind; 
however, the discussion in this section is not limited to that
particular instance. 

\emph{Projection methods} solve \eqref{e:cfp} 
by generating a sequence of vectors by
using the projectors $P_{C_1},\ldots,P_{C_m}$ onto the individual sets.
(For further background material on projection methods,
see, e.g., \cite{BC2011}, \cite{Cegielski}, and \cite{CZ}.)
As is illustrated by Section~\ref{sec:CPO}, these projectors are
available for a variety of constraints appearing in practical
applications.
Before we proceed to catalogue the algorithms, we note that 
some algorithms (see, e.g., Algorithm~\ref{a:D--R} below) 
work intrinsically only
with \emph{two} constraint sets. 
Because $m$ in our application is not large, this turns out to be 
not a handicap at all as we can 
reformulate \eqref{e:cfp} in a product space in the following fashion: 
In the \emph{Hilbert product space}
\begin{equation}
\bX := X^m,
\end{equation}
equipped with
\begin{equation}
\scal{\bx}{\by} = \sum_{i=1}^{m}\scal{x_i}{y_i}
\quad\text{and}\quad
\|\bx\| = \sqrt{\sum_{i=1}^m \|x_i\|^2},
\end{equation}
where $\bx = (x_1,\ldots,x_m)$ and $\by=(y_1,\ldots,y_m)$ belong to
$\bX$, 
we consider the \emph{Cartesian product} of the constraints together 
with the \emph{diagonal} in $\bX$,
i.e., 
\begin{equation}
\bC := C_1\times\cdots\times C_m
\quad\text{and}\quad
\bD := \menge{(x,\ldots,x)\in \bX}{x\in X}.
\end{equation}
Then for every $x\in X$, we have the key equivalence 
\begin{equation}
x\in C 
\quad\Leftrightarrow\quad
(x,\ldots,x)\in\bC\cap\bD;
\end{equation}
thus, the original $m$-set feasibility problem \eqref{e:cfp} in $X$ 
reduces to \emph{two-set} feasibility problem for the sets 
$\bC$ and $\bD$ in $\bX$. 
This approach is viable because the projectors onto $\bC$ and
$\bD$ are given explicitly 
(see, e.g., \cite[Propositions 28.3 and 28.13]{BC2011}) by
\begin{subequations}
\label{e:Pproduct}
\begin{equation}
P_\bC \colon (x_1,\ldots,x_m)\mapsto (P_{C_1}x_1,\ldots,P_{C_m}x_m)
\end{equation}
and by
\begin{equation}
P_\bD \colon (x_1,\ldots,x_m)\mapsto (y,\ldots,y),
\quad\text{where $y=\tfrac{1}{m}(x_1+\cdots+x_m)$,}
\end{equation}
\end{subequations}
respectively. 

\subsection{Projection methods and\ \ldots\ swiss army knives!}
\label{sec:swissarmy}

Projection methods use the projectors onto the given constraint sets
in some fashion. Because the squared distance functions,
$d^2(\cdot,C_i)$, are Fr\'echet differentiable with derivative
$2(\Id-P_{C_i})$ (see, e.g., \cite[Corollary~12.30]{BC2011}) projection
methods are \emph{first-order methods} --- these are known
to allow at best for linear convergence. 
Not surprisingly, they are not always competitive with 
special-purpose solvers \cite{Gould}; 
however, when projection methods succeed 
(see \cite{CCCDH} for a compelling set of examples), 
then they have a variety of \emph{very attractive features}:

\begin{itemize}
\item Projection methods are \emph{easy to understand.}
This is important an industry, where mathematical/algorithmic
considerations are only one part of an engineer's job. The engineer
will not typically be familiar with the latest research developments
in all branches of relevant mathematics.
On the other hand, the idea of a projection
method is often immediately grasped by drawing some sketches. 
\item Projection methods are \emph{easy to implement.}
In Section~\ref{sec:CPO}, we have seen various formulas for projection
methods. These involve simple operations from linear algebra and are
easily programmed.
\item Projection methods are \emph{easy to maintain.}
Once implemented, the code for these methods is typically small and in
fact smaller than other pieces of code dealing with data input/output.
This makes maintenance straightforward. 
\item Projection methods are \emph{easy to deploy.}
Because of typically small memory requirements, it makes them much
easier to deploy on low memory computers like mobile devices.  
Also, the code base required for projection methods is in most cases
significantly smaller than the size of libraries for linear or
nonlinear optimization solver software.
Thus, projection methods satisfy some of the 
key requirements for \emph{embedded optimization} \cite{Boyd2012}, 
where the solution of one method is used within an encompassing algorithm.
\item Projection methods are \emph{inexpensive to implement.}
Because of typically straight forward implementations, there is no
need for commercial optimization solver software.
\item Projection methods can be \emph{very fast}. 
If the iterations in a projection method can be executed
quickly, then for certain classes of problems 
projection methods can become very competitive with traditional optimization
algorithms.  
In Section \ref{sec:lpsimplex} below, we illustrate 
the enormous potential of projection methods when compared to 
algorithms for linear programming or even 
mixed-integer linear programming.
\end{itemize}

In summary, \emph{projection methods possess the same essential 
characteristics of swiss army knives}: 
they are flexible, lightweight, simple and very convenient provided
they get the job done.
If the saw included with the swiss army knife in your pocket 
cuts the branch of the tree, then there is no need for you to either 
get a big saw out of the garage or to buy a chainsaw from the hardware
store! 
The road design feasibility problem analyzed in this paper adds
a new compelling success story.

\subsection{A catalogue of projection methods}

In this subsection, we provide a list of projection methods.
Each of these methods produces a sequence that 
converges --- sometimes after applying a suitable operator ---
to a point in $C$ provided that $C\neq\varnothing$
(and perhaps an additional assumption is satisfied). 
Basic references are the books
\cite{BC2011}, \cite{Cegielski}, \cite{CZ}, \cite{Deutsch},
\cite{GK}, and \cite{GR}; if these books are insufficient, we
include further pointers in accompanying remarks. 

Note that numerous generalizations of projection 
methods are known. These typically involve 
additional parameters (e.g., weights and 
relaxation parameters).
Because we shall compare these methods numerically, 
we have to restrict our attention to the most basic 
instance of each method.

While these methods are generally not able to find $P_C\anfang$,
they may lead to feasible solutions that are ``fairly close'' to
$\anfang$ provided the starting point is chosen to be $\anfang$.

We start with the method of cyclic projections, which has a long
history (see, e.g., \cite{Deutsch}).

\begin{algo}[cyclic projections (CycP)]
\label{a:CycP}
Set $x_0=\anfang$ and update
\begin{equation}
\label{e:CycP}
(\forall k\in\NN)\quad
x_{k+1} := Tx_{k},
\quad\text{where}\quad
T := P_{C_m}P_{C_{m-1}}\cdots P_{C_1}. 
\end{equation}
\end{algo}

A modern variant of CycP replaces the projectors
in \eqref{e:CycP} by intrepid counterparts when available. 

\begin{algo}[cyclic intrepid projections (CycP+)]
\label{a:CycP+}
Set $x_0=\anfang$ and update
\begin{equation}
\label{e:CycP+}
(\forall k\in\NN)\quad
x_{k+1} := Tx_{k},
\quad\text{where}\quad
T := R_{m}R_{m-1}\cdots R_{1}, 
\end{equation}
where each $R_i$ is $P_{C_i}$ or an intrepid counterpart thereof. 
\end{algo}

\begin{remark}[CycP+]
CycP+ (also known as ART3+) is available for \eqref{e:www} because we have
the intrepid counterparts of the projectors at our disposal
(see Propositions~\ref{p:iPS} and \ref{p:iPcurv}, and
Footnote~\ref{fn:PY}). 
This method is fairly recent; see
\cite{CCCDH}, \cite{CCP}, and \cite{HC} (for constructions based on
halfspaces and subgradient projectors), and also
\cite{BK} (where one of the sets is an \emph{obtuse cone} $K$ 
so that the intrepid projector is actually the reflector $2P_K-\Id$). 
\end{remark}

While Algorithm~\ref{a:CycP} is \emph{sequential}, 
the following method is \emph{parallel}:

\begin{algo}[parallel projections (ParP)]
\label{a:ParP}
Set $x_0=\anfang$ and update
\begin{equation}
\label{e:ParP}
(\forall k\in\NN)\quad 
x_{k+1} := Tx_{k},
\quad\text{where}\quad
T := \tfrac{1}{m}\big(P_{C_1}+\cdots + P_{C_m}\big). 
\end{equation}
\end{algo}

\begin{remark}[ParP]
In view of \eqref{e:Pproduct},
it is interesting to note that 
ParP is equivalent to iterating $P_{\bD}P_{\bC}$, i.e.,
to applying CycP to the subsets $\bC$ and $\bD$ of $\bX$. 
See also \cite[Corollary~2.6]{Reich83} and 
\cite[Section~6]{BBDyk}. 
\end{remark}

The next method can be seen as a hybrid of CycP and ParP.

\begin{algo}[string-averaging projections (SaP)]
Set $x_0 = \anfang$ and update
\begin{equation}
\label{e:SaP}
(\forall k\in\NN)\quad 
x_{k+1} := Tx_{k},
\quad\text{where}\quad
T := \tfrac{1}{m}\big(P_{C_1}+P_{C_2}P_{C_1}+\cdots + P_{C_m}\cdots
P_{C_2}P_{C_1}\big). 
\end{equation}
\end{algo}

\begin{remark}[SaP]
For further information, see 
\cite{CEH} and \cite{CT} 
(and also \cite[Example~2.14]{BB96}).  
\end{remark}

\begin{algo}[extrapolated parallel projections (ExParP)]
Set $x_0=v$ and update, for every $k\in\NN$, 
\begin{equation}
\label{e:ExParP}
x_{k+1} := Tx_{k}, 
\quad\text{where}\quad 
Tx := \begin{cases}
\displaystyle 
x +
\frac{\sum_{i=1}^m\|x-P_{C_i}x\|^2}{\|\sum_{i=1}^m(x-P_{C_i}x)\|^2}\sum_{i=1}^m(P_{C_i}x-x),
&\text{if $x\notin C$;}\\
x, &\text{otherwise.}
\end{cases}
\end{equation}
\end{algo}

\begin{remark}[ExParP]
This method is actually an instance of the \emph{subgradient projection
method} applied to the function $x\mapsto \sum_{i=1}^{m}d^2(x,C_i)$;
see \cite{C97} for further information. 
\end{remark}

\begin{algo}[extrapolated alternating projections (ExAltP)]
Assume that $C_1$ is an affine subspace --- this is the case in
our application when we choose the interpolation constraint. 
Set $x_0=\anfang$, and let $k\in\NN$. 
Let $I_k$ be a nonempty subset of $\{1,\ldots,m\}$ 
containing exactly $m_k$ indices such that 
for each $j\in\{2,\ldots,m\}$, 
$j$ belongs to $I_k$ frequently. 
Given $x_k$, 
set 
\begin{subequations}
\label{e:ExAltP}
\begin{equation}
z_k = P_{C_1}x_k,
\;\;
p_k = P_{C_1}\bigg(\frac{1}{m_k}\sum_{i\in I_k}P_{C_i}z_k\bigg),
\;\;
\mu_k = 
\begin{cases}
\displaystyle 
\frac{\sum_{i\in I_k}\|z_k-P_{C_i}z_k\|^2}{m_k\|p_k-z_k\|^2},
&\text{if $\displaystyle z_k\notin\bigcap_{i\in I_k}C_i$};\\
1, &\text{otherwise,}
\end{cases}
\end{equation}
and then update
\begin{equation}
x_{k+1} = z_k + \mu_k(p_k-z_k).
\end{equation}
\end{subequations}
\end{algo}

\begin{remark}[ExAltP]
See \cite[Algorithm~3.5]{BCK} for further information on this method.
In our implementation, we chose 
$I_k=\{2,\ldots,m\}$ so that  $m_k=m-1$. 
\end{remark}

The next method is different from the previous ones
in two aspects: 
First, it truly operates in $\bX$ and thus
has increased storage requirements --- fortunately, $m$ is small in our
application so this is of no concern when dealing with \eqref{e:www}.
Second, the sequence of interest is actually different and derived
from another sequence that governs the iteration. 

\begin{algo}[Douglas--Rachford (D--R)]
\label{a:D--R}
Set $\bx_0=(\anfang,\ldots,\anfang)\in\bX=X^m$.
Given $k\in\NN$ and
$\bx_k = (x_{k,1},\ldots,x_{k,m})\in \bX$, update to 
$\bx_{k+1} = (x_{k+1,1},\ldots,x_{k+1,m})$, where
\begin{subequations}
\label{e:D--R} 
\begin{equation}
\label{e:D--Rmon}
\bar{x}_k = \frac{1}{m}\sum_{i=1}^{m}x_{k,i}
\end{equation}
and 
\begin{equation}
(\forall i\in\{1,\ldots,m\})\;\;
x_{k+1,i} = x_{k,i} - \bar{x}_k + P_{C_i}(2\bar{x}_k-x_{k,i}).
\end{equation}
\end{subequations}
The sequence of interest is not $(\bx_k)_{k\in\NN}$ but rather
$(\bar{x}_k)_{k\in\NN}$. 
\end{algo}

\begin{remark}[general Douglas--Rachford algorithm]
\label{r:genD--R}
Let us briefly sketch how the update formula \eqref{e:D--R} 
is derived from the general splitting version of
D--R, which aims to minimize the sum $f+g$ of proper lower
semicontinuous convex functions $f\colon X\to\RX$
and $g\colon X\to\RX$. 
Given $x_0\in X$, the algorithm proceeds via
\begin{equation}
\label{e:genD--R}
(\forall\knn)\quad
x_{k+1} = Tx_k,\;\;
\text{where}\;\;
T = \prox_{g}(2\prox_{f}-\Id)+\Id-\prox_{f},
\end{equation}
and where $\prox_{f}$ denotes the \emph{proximal mapping} (or
\emph{proximity operator}) of $f$, i.e., 
$\prox_{f}(y)$ is the unique minimizer of the function
$x\mapsto \thalb\|x-y\|^2 + f(x)$; the sequence to monitor is
$(\prox_{f}x_k)_\knn$.
Now assume that $A$ and $B$ are two nonempty closed convex subsets of
$X$. 
The \emph{indicator function} $\iota_A$ of $A$ takes the value
$0$ on $A$, and $+\infty$ outside $A$, and analogously for $\iota_B$. 
We set $f=\iota_B$ and $g=\iota_A$.  It is then clear that 
\begin{equation}
\prox_{f} = P_B 
\quad\text{and}\quad
\prox_{g} = P_A 
\end{equation}
and that \eqref{e:genD--R} turns into
\begin{equation}
(\forall\knn)\quad
x_{k+1} = Tx_k,\;\;
\text{where}\;\;
T = P_A(2P_B-\Id)+\Id-P_B
\end{equation}
Applying this in $\bX$ (with $A=\bC$ and $B=\bD$) and
recalling \eqref{e:Pproduct}, 
we obtain \eqref{e:D--R}. 
Viewed directly in $\bX$, 
we obtain the iteration 
\begin{subequations}
\begin{equation}
\label{e:bD--R} 
\bx_{k+1} = \bT\bx_{k},
\;\;\text{where}\;\;
\bT := P_{\bC}(2P_{\bD}-\bId)+\bId - P_{\bD}
= \frac{\bId + (2P_\bC-\bId)(2P_\bD-\bId)}{2},
\end{equation}
and we monitor the sequence
\begin{equation}
(P_{\bD}\bx_k)_{k\in\NN}.
\end{equation}
\end{subequations}
For convergence results, see 
\cite[Chapter~26]{BC2011},
\cite{EckBer}, and
\cite{LM}. 
\end{remark}

\begin{remark}[D--R vs alternating direction method of multipliers
(ADMM)] \ 
\label{r:DRvsADMM} \\
ADMM is a very popular method \cite{Boyd} that can also be adapted
to solve feasibility problems for two sets.
Suppose we wish to find a point in $A\cap B$, where
$A$ and $B$ are two nonempty closed convex subsets of $X$.
Given $u_0\in X$ and $b_0\in X$, ADMM generates three sequences
$(a_k)_{k\geq 1}$, $(b_k)_\knn$, and $(u_k)_\knn$ via
\begin{equation}
(\forall\knn)\quad
a_{k+1}:=P_A(b_k-u_k),
\;
b_{k+1}:=P_B(a_{k+1}+u_k),
\;
u_{k+1}:=u_k+a_{k+1}-b_{k+1}. 
\end{equation}
On the other hand, D--R for this problem, 
with starting point $x_0\in X$,  produces two 
sequences $(x_k)_\knn$ and $(y_k)_\knn$ via
\begin{equation}
(\forall\knn)\quad y_k := P_Bx_k,\;x_{k+1}:=P_A(2y_k-x_k)+x_k-y_k.
\end{equation}
Now assume that
\begin{equation}
x_0 = b_0 \in B
\;\text{and}\; 
u_0 = 0.
\end{equation}
Then 
$y_0 = P_Bx_0 = x_0$,
$x_1 = P_A(2y_0-x_0)+x_0-y_0 = P_Ax_0
=P_A(b_0-0)=P_A(b_0-u_0)=a_1 = a_1+u_0$
and
$y_1 = P_Bx_1 = P_B(a_1+u_0)=b_1$.
Furthermore, assume that for some $k\geq 1$, we have 
$x_k=a_k + u_{k-1}$ and $y_k = b_k$.
Then
$2y_k-x_k=2b_k-(u_{k-1}+a_k)=b_k-a_k+b_k-u_{k-1}
=u_{k-1}-u_k+b_k-u_{k-1}=b_k-u_k$
and $x_k-y_k=u_{k-1}+a_k-b_k=u_k$;
in turn, this implies 
$x_{k+1}=P_A(2y_k-x_k)+x_k-y_k=P_A(b_k-u_k)+u_k
=a_{k+1}+u_k$
and $y_{k+1}=P_Bx_{k+1}=P_B(a_{k+1}+u_k)=b_{k+1}$. 
It follows inductively that\footnote{\ 
In fact, this argument works
much more generally when the projectors $P_A$ and $P_B$ are replaced by 
arbitrary \emph{resolvents} and $b_0$ is a fixed point of the second
resolvent.}
\begin{equation}
(x_k)_{k\geq 1} = (a_k+u_{k-1})_{k\geq 1}
\;\;\text{and}\;\;
(y_k)_{k\geq 1} = (b_k)_{k\geq 1}.
\end{equation}
In this sense, 
\begin{equation}
\text{
D--R and ADMM \emph{are equivalent}. 
}
\end{equation}
See also \cite{Boyd,CP,GM,GlM} for further information on ADMM and
related methods. 
Furthermore, if $A$ and $B$ are linear subspaces, then $(x_k)_\knn$,
the sequence governing DR, has the update rule
\begin{equation}
x_{k+1} = (P_AP_B+P_{A^\perp}P_{B^\perp})x_k
\end{equation}
since $P_A(2P_B-\Id)+\Id-P_B = P_A(2P_B-P_B-P_{B^\perp}) + P_{B^\perp}
=P_A(P_B-P_{B^\perp}) + P_{A}P_{B^\perp} + P_{A^\perp}P_{B^\perp} =
P_AP_B + P_{A^\perp}P_{B^\perp}$.
Using, e.g., \cite[Corollary~3.9]{BCL}, one  can further show
that $\Fix (P_AP_B+P_{A^\perp}P_{B^\perp}) = (A\cap
B)+(A^\perp\cap B^\perp)$. 
\end{remark}

\subsection{Summary of feasibility algorithms}

\label{ss:summaryfeas}

\begin{table*}[H] \centering
\begin{tabular}{@{}lllr@{}} \toprule
Name &Acronym & Formula & Monitor\\ \midrule
Cyclic Projections  &CycP& \eqref{e:CycP} & $(x_k)_{k\in\NN}$ \\
Cyclic Intrepid Projections  &CycP+& \eqref{e:CycP+} & $(x_k)_{k\in\NN}$ \\
Parallel Projections &ParP      & \eqref{e:ParP} & $(x_k)_{k\in\NN}$\\
String-averaging Projections &SaP      & \eqref{e:SaP} & $(x_k)_{k\in\NN}$\\
Extrapolated Parallel Projections & ExParP & \eqref{e:ExParP} & $(x_k)_{k\in\NN}$  \\
Extrapolated Alternating Projections & ExAltP & \eqref{e:ExAltP} & $(x_k)_{k\in\NN}$  \\
Douglas--Rachford & D--R & \eqref{e:D--R} & $(\bar{x}_k)_{k\in\NN}$ \\ 
\bottomrule
 \end{tabular}
\end{table*}

Note that all algorithms in this table proceed by iterating an
operator, and monitoring either the iterates directly, or some simple
version thereof. This is a key point when we revisit these algorithms
in the next section.

\subsection{Superiorization: between feasibility and best approximation}

\label{ss:superior}

The algorithms considered so far 
are designed to solve 
the \emph{feasibility problem} \eqref{e:cfp}.
Let $\anfang\in X\smallsetminus C$. 
We shall discuss algorithms for 
finding $P_C\anfang$ in Section~\ref{sec:ba} below.
This \emph{best approximation problem} 
is equivalent to\footnote{\ We work here with a
squared version of \eqref{e:distance} because the objective function
is then differentiable.}
solving the optimization problem
\begin{equation}
\label{e:bap}
\min_{x\in C} \thalb\|x-\anfang\|^2.
\end{equation}
The new paradigm of \emph{superiorization} (see, e.g., \cite{CDH}) lies
between feasibility and this best approximation problem. 
It is not quite trying to solve \eqref{e:bap}; rather, the objective
is to find a feasible point that is a superior to one returned by a
feasibility algorithm. To explain this in detail, we
assume that $T\colon X\to X$ satisfies 
\begin{equation}
\Fix T = C;
\end{equation}
hence, \eqref{e:bap} is equivalent to 
\begin{equation}
\label{e:bap'}
\min_{x\in \Fix T} \thalb\|x-\anfang\|^2.
\end{equation}
Applying the superiorization approach to \eqref{e:bap'},
we obtain the following abstract algorithm: 

\begin{equation}
\begin{minipage}{.65\textwidth}
\begin{algorithm}[H]
\TitleOfAlgo{Superiorization of $T$}
\SetAlgoLined
\DontPrintSemicolon
\KwData{$\anfang\in X,\varepsilon>0$}
\KwResult{$x_k$}
$k \leftarrow 0$\;
$x_0 \leftarrow \anfang$\;
$\theta \leftarrow 1$\;
\While{$d(x_k) > \varepsilon$}{
\eIf{$\|x_k-\anfang\| > 0$}
{$\tilde{x} \leftarrow x_k + \theta(x_k-\anfang)/\|x_k-\anfang\|$\;}
{$\tilde{x} \leftarrow x_k$\;}
$\theta \leftarrow \theta/2$\;
\If{$\|\tilde{x}-\anfang\|\leq\|x_k-\anfang\|$ \KwSty{and}~ $d(T\tilde{x}) < d(x_k)$}
{$x_{k+1} \leftarrow T\tilde{x}$\;}
$k \leftarrow k+1$\;
}
\end{algorithm}
\end{minipage}
\end{equation}
Note that $d\colon X\to\RR_+$ is a performance function
satisfying $d(x)=0$ if and only if $x\in C$.
(In Section~\ref{sec:numeresul}, we use \eqref{e:d}.)
With the exception of D--R, each algorithm in
Section~\ref{ss:summaryfeas} has a \emph{superiorized} counterpart. 
(It is not clear  how D--R should be
superiorized because the fixed point set of the operator $T$
governing the iteration is generally different from $C$, the set of interest.)
We denote the superiorized counterpart of CycP by sCycP, and
analogously for the other algorithms. 
For remarks on the numerical performance of these algorithms, see 
Section~\ref{ss:ressuper}. 

\section{Best approximation algorithms}
\label{sec:ba}

\subsection{The problem and notation}

We continue to assume that 
$m\in\{2,3,\ldots\}$ and that
$C_1,\ldots,C_m$ are closed convex subsets of $X$ such that 
\begin{equation}
C := C_1\cap\cdots\cap C_m\neq\varnothing.
\end{equation}
Let $\anfang\in X$. 
We wish to determine
\begin{equation}
P_{C}\anfang.
\end{equation}
Before we present a selection of pertinent best approximation
algorithms, let us fix some notation for ease of use. 
It is notationally convenient to 
occasionally work with cyclic remainders
\begin{equation}
\label{e:cycrem}
[1]_m = 1,
[2]_m = 2,\ldots,[m]_m=m,[m+1]_m=1, \ldots,
\end{equation}
taken from $\{1,\ldots,m\}$. 
We also require the operator $Q$ defined by 
\footnote{\ If $\rho=0$ and $\chi<0$, then the output of $Q$ 
is undefined
--- this corresponds to the case when $C=\varnothing$.}
\begin{align}
\label{e:Q}
Q\colon X\times X\times X &\to X\\
(x,y,z) &\mapsto 
\begin{cases}
z, &\text{if $\rho=0$ and $\chi\geq 0$;}\\
x + \Big(1+\frac{\chi}{\nu}\Big)(z-y),
&\text{if $\rho > 0$ and $\chi\nu\geq \rho$;}\\
y + \frac{\nu}{\rho}\Big(\chi(x-y)+\mu(z-y)\Big),
&\text{if $\rho>0$ and $\chi\nu<\rho$,}
\end{cases} \nonumber
\end{align}
where 
$\chi := \scal{x-y}{y-z}$,
$\mu := \|x-y\|^2$,
$\nu := \|y-z\|^2$,
and $\rho := \mu\nu-\chi^2$. 
An analogous formula holds for $\bQ\colon\bX\times\bX\times\bX\to\bX$.

\subsection{A catalogue of best approximation methods}

As in Section~\ref{sec:FPM}, we present a list of best
approximation methods based on projectors and comments. 
Unless stated otherwise, 
each of these methods produces a main/governing sequence that is
converging to $P_C\anfang$.

\begin{algo}[Halpern--Wittmann (H--W)]
\label{a:H-W}
Set $x_0=\anfang$ and update
\begin{equation}
\label{e:H-W}
(\forall\knn)\quad 
x_{k+1} := \tfrac{1}{k+1}\anfang + \tfrac{k}{k+1}
P_{C_m}P_{C_{m-1}}\cdots P_{C_1}x_k. 
\end{equation}
\end{algo}

\begin{remark}[H--W]
This algorithm was introduced by Halpern \cite{Halpern} while
Wittmann \cite{Wittmann} proved convergence for the choice of parameters
presented in Algorithm~\ref{a:H-W}. 
Many variants have been proposed and studied. 
\end{remark}

\begin{algo}[Cyclic Dykstra algorithm (CycDyk)]
Set $x_0 := \anfang$, 
$q_{-(m-1)} := q_{-(m-2)} := \cdots := q_{-1} := q_0 := 0$,
and update
\begin{equation}
\label{e:CycDyk}
(\forall k\in\NN)\quad
x_{k+1} := P_{C_{[k+1]_m}}(x_k + q_{k+1-m})
\;\;\text{and}\;\;
q_{k+1} := x_{k} + q_{k+1-m} - x_{k+1}.
\end{equation}
\end{algo}

\begin{remark}[CycDyk]
For convergence proofs, see, e.g., 
\cite{B-D} or \cite[Theorem~29.2]{BC2011}.
See also \cite[Section~3]{CDV} for connections to the
\emph{forward--backward method}. 
\end{remark} 

The following method operates in $\bX$.

\begin{algo}[Parallel Dykstra algorithm (ParDyk)]
\label{a:ParDyk}
Set
$(\by_0,\bz_0)=(\anfang,\ldots,\anfang,0,\ldots,0)\in \bX\times\bX$, 
and let
$(\by_k,\bz_k)=(y_{k,1},\ldots,y_{k,m},z_{k,1},\ldots,z_{k,m})\in
\bX\times\bX$ be given. Then the next iterate is
$(\by_{k+1},\bz_{k+1})=(y_{k+1,1},\ldots,y_{k+1,m},z_{k+1,1},\ldots,z_{k+1,m})$,
where the sequence $(\bar{x}_k)_{k\in\NN}$ to monitor is 
\begin{subequations}
\label{e:ParDyk}
\begin{equation}
\bar{x}_k = \frac{1}{m}\sum_{i=1}^{m}y_{k,i}
\end{equation}
and the update formulas are 
\begin{equation}
\label{e:ParDykii}
(\forall i\in\{1,\ldots,m\})\quad
y_{k+1,i} :=  P_{C_i}(z_{k,i}+\bar{x}_k)
\;\;\text{and}\;\;
z_{k+1,i} := z_{k,i} + \bar{x}_k - y_{k+1,i}.
\end{equation}
\end{subequations}
\end{algo}

\begin{remark}[ParDyk]
ParDyk is CycDyk applied to the subsets $\bC$ and $\bD$ of $\bX$; 
see, e.g., \cite[Theorem~6.1]{BBDyk}. 
\end{remark}

\begin{remark}[Dykstra vs ADMM]
\label{r:DykvsADMM}
Let $A$ and $B$ be two nonempty closed convex subsets of $X$,
and let $\anfang\in X$.
CycDyk for finding $P_{A\cap B}w$ generates
sequences $(a_k)_{k\geq 1}$, $(b_{k})_\knn$,
$(p_k)_\knn$, and $(q_k)_\knn$ as follows\footnote{\ Here the sequence
$(q_k)_\knn$ from \eqref{e:CycDyk} is split into subsequences
corresponding to odd and even terms for easier readability.}:
Set $b_0 := \anfang$, $p_0 := q_0 :=0$, and
for every $\knn$, update
\begin{subequations}
\label{e:BieberDyk}
\begin{align}
a_{k+1} &:= P_A(b_k+p_k),\;\;  & p_{k+1} &:= b_k+p_k-a_{k+1},\;\;\\
b_{k+1} &:= P_B(a_{k+1}+q_k),\;\; & q_{k+1} &:= a_{k+1}+q_{k}-b_{k+1}.
\end{align}
\end{subequations}
On the other hand,
given $u_0\in X$ and $y_0\in X$, ADMM for finding a point in $A\cap B$
--- not necessarily $P_{A\cap B}w$ --- generates three sequences
$(x_k)_{k\geq 1}$, $(y_k)_\knn$, and $(u_k)_\knn$ via
\begin{equation}
\label{e:BieberADMM}
(\forall\knn)\quad
x_{k+1}:=P_A(y_k-u_k),
\;
y_{k+1}:=P_B(x_{k+1}+u_k),
\;
u_{k+1}:=u_k+x_{k+1}-y_{k+1}. 
\end{equation}
Now let us assume that 
\begin{equation}
\text{$B$ is a linear subspace}
\end{equation}
and that $u_0\in B^\perp$. 
Then it is well known that the Dykstra update \eqref{e:BieberDyk} simplifies
to 
\begin{equation}
\label{e:BieberDyklin}
a_{k+1} := P_A(b_k+p_k),\;\;
b_{k+1} := P_B(a_{k+1}),\;\;
p_{k+1} := p_k+b_k-a_{k+1},
\end{equation}
because the sequence $(q_k)_\knn$ lies in $B^\perp$ and thus becomes
``invisible'' when computing $(b_k)_\knn$ due to the linearity of
$P_B$. Turning to ADMM, we observe that $(u_k)_\knn$ lies in
$B^\perp$ which simplifies the update for $y_{k+1}$ to
$y_{k+1} := P_B(x_{k+1})$. Setting $(v_k)_\knn := -(u_k)_\knn$, 
we see that \eqref{e:BieberADMM} turns into 
\begin{equation}
\label{e:BieberADMMlin}
(\forall\knn)\quad
x_{k+1}:=P_A(y_k+v_k),
\;
y_{k+1}:=P_B(x_{k+1}),
\;
v_{k+1}:=v_k+y_{k+1}-x_{k+1}. 
\end{equation}
Comparing \eqref{e:BieberDyklin} to \eqref{e:BieberADMMlin},
we see that the update formulas look \emph{almost} identical:
The only difference lies in the update formulas of the auxiliary sequences
$(p_k)_\knn$ and $(v_k)_\knn$ --- the former works with $b_{k}$ while
the latter incorporates immediately the more recent update $y_{k+1}$.
However, the resulting sequences and hence algorithms appear to be 
\emph{different}\footnote{\label{fn:goof} 
It is stated on \cite[page~34f]{Boyd}
that (in our notation) ``\eqref{e:BieberADMMlin} is exactly Dykstra's
alternating projections method \ldots which is far more
efficient than the classical method [of alternating projections] that
does not use the dual variable $v$.''. 
This statement appears to be at least ambiguous because the crucial
starting points are not specified.}.
Let us now further specialize by additionally assuming that 
\begin{equation}
\text{$A$ is a linear subspace.}
\end{equation}
Then the Dykstra update \eqref{e:BieberDyklin} does not require the
sequence $(p_k)_\knn$ anymore (because it lies in $A^\perp$ and it
thus plays no role in the generating of $(a_k)_{k\geq 1}$), 
and it further simplifies to
\begin{equation}
\label{e:wewin5}
a_{k+1} := P_A(b_k),\;\;
b_{k+1} := P_B(a_{k+1}).
\end{equation}
Hence, as is well known, 
Dykstra turns into CycP, which is also 
known as \emph{von Neumann's alternating
projection method} in this special case. 
On the other hand, these additional assumptions do not seem to
simplify \eqref{e:BieberADMMlin}. In fact, the behaviour of Dykstra
can starkly differ from ADMM even in this setting:
Suppose that $X=\RR^2$, that $A=\RR\cdot(1,1)$ and $B =
\RR\times\{0\}$. Then the sequence $(b_k)_\knn$ with starting point
$b_0 := v:=(1,0)$ turns out to be $(2^{-k},0)_{\knn}$. 
In contrast (see also Footnote~\ref{fn:goof}), 
we compute the following ADMM updates, where
$y_0 := v:= (1,0)$ and $v_0 := (0,0)$:
$x_1 = P_A(y_0+v_0) = P_A(1,0)=\thalb(1,1)$,
$y_1 = P_B(x_1)= \thalb(1,0)$,
$v_1 = v_0 +y_1-x_1 = \thalb(0,-1)$,
$x_2 = P_A(y_1+v_1) = P_A(\thalb(1,-1)) = (0,0)$,
and $y_2 = P_B(x_2) = (0,0)$. 
Since $x_2 = y_2\in A\cap B$, the algorithm terminates whenever
feasibility is the implemented stopping criterion. 
\end{remark}

\begin{algo}[Haugazeau's algorithm with cyclic projections (hCycP)]
Set $x_0=\anfang$, and update\footnote{\ Recall \eqref{e:cycrem}.}
\begin{equation}
\label{e:Haug}
(\forall\knn)\quad 
x_{k+1} := Q(x_0,x_k,P_{C_{[k+1]_m}}x_k).
\end{equation}
\end{algo}
\begin{remark}[hCycP]
For convergence proofs, see \cite{Haug} or
\cite[Corollary~29.8]{BC2011}.
The general pattern for methods that undergo a Haugazeau-type
modification is  (see \cite{MOR} for details)
\begin{equation}
(\forall\knn)\quad 
x_{k+1} := Q(x_0,x_k,T_{\mathrm{i}(k)}x_k),
\end{equation}
where $(T_i)_{i\in I}$ is a family of operators and 
$\mathrm{i}\colon\NN\to I$ selects which operator is drawn at
iteration $k$. This gives rise to
many variants; in the following, we shall focus on a representative
selection.
\end{remark}

Here is a Haugazeau-type modification of ParP
(Algorithm~\ref{a:ParP}): 

\begin{algo}[hParP]
Set $x_0:=\anfang$ and update 
\begin{equation}
\label{e:hParP}
(\forall\knn)\quad 
x_{k+1} := Q(x_0,x_k,Tx_k),
\quad
\text{where}\;\;
T = \tfrac{1}{m}\big(P_{C_1}+\cdots + P_{C_m}\big). 
\end{equation}
\end{algo}

The following Haugazeau-type modification of 
D--R operates in the product space $\bX$.

\begin{algo}[hD--R]
Let $\bT\colon\bX\to\bX$ be the operator governing D--R 
(see \eqref{e:bD--R}), and set $\bx_0=(\anfang,\ldots,\anfang)\in\bX$.
Given $\knn$ and $\bx_k=(x_{k,1},\ldots,x_{k,m}) \in\bX$, update 
the governing sequence by 
\begin{subequations}
\label{e:hD--R}
\begin{equation}
\bx_{k+1} := \bQ(\bx_0,\bx_k,\bT\bx_k),
\end{equation}
and monitor 
\begin{equation}
\bar{x}_k := \frac{1}{m}\sum_{i=1}^{m}x_{k,i}.
\end{equation}
\end{subequations}
\end{algo}

\begin{remark}[hD--R]
To show that $\bar{x}_k\to P_C(\anfang)$, 
use \cite[Example~8.2]{BBHM}. 
\end{remark}

The next algorithm is a variant of D--R tailored for best
approximation. It also operators in the product space $\bX$.

\begin{algo}[baD--R]
\label{a:D--Rba}
Set $\bx_0 = (\anfang,\ldots,\anfang)\in\bX$.
Given $\knn$ and 
$\bx_k = (x_{k,1},\ldots,x_{k,m})\in \bX$, 
update to $\bx_{k+1} = (x_{k+1,1},\ldots,x_{k+1,m})$, where
\begin{subequations}
\label{e:D--Rba} 
\begin{equation}
\label{e:D--Rbamon}
\bar{x}_k = \frac{1}{m}\sum_{i=1}^{m}x_{k,i}
\end{equation}
and 
\begin{equation}
(\forall i\in\{1,\ldots,m\})\;\;
x_{k+1,i} = x_{k,i} - \bar{x}_k +
P_{C_i}\big((\anfang+2\bar{x}_k-x_{k,i})/2\big).
\end{equation}
\end{subequations}
\end{algo}

\subsection{Remarks on baD--R}
\label{ss:baD--R}

We comment on various aspects of baD--R and start with its genesis. 

\begin{remark}[baD--R]
Let $A$ and $B$ be nonempty closed convex subsets of $X$,
and let $\anfang\in X$. 
Revisit Remark~\ref{r:genD--R}, in which we showed how the
feasibility version of D--R, Algorithm~\ref{a:D--R}, arises from the
general problem of minimizing $f+g$, by choosing $f=\iota_B$ and
$g=\iota_A$. 
To explain Algorithm~\ref{a:D--Rba}, we keep $f = \iota_B$ for which
\begin{equation}
\label{e:wewin1}
\prox_{f} = P_B.
\end{equation}
However, this time we take 
$g\colon x\mapsto \thalb\|x-\anfang\|^2 + \iota_A(x)$.
Then the (unique) minimizer of $f+g$ is indeed $P_{A\cap B}\anfang$.
Let $x$ and $y$ be in $X$. Then $y= \prox_{g}x$
$\Leftrightarrow$
$x\in(\Id+\partial g)(y) = 2y-\anfang+N_A(y)$
$\Leftrightarrow$
$\thalb(x+\anfang)\in y + \thalb N_A(y) = y+N_A(y)$ and so 
$y=P_A(\thalb(x+\anfang))$; that is,
\begin{equation}
\label{e:wewin2}
\prox_g\colon x\mapsto P_A\big(\thalb(x+\anfang)\big).
\end{equation}
In view of \eqref{e:wewin1} and \eqref{e:wewin2}, 
the general update formula for the Douglas--Rachford splitting method, 
\eqref{e:genD--R}, becomes 
\begin{equation}
\label{e:wewin3}
(\forall\knn)\quad
x_{k+1} = Tx_k,\;\;
\text{where}\;\;
T \colon x\mapsto  P_A((2P_Bx-x+v)/2)+ x-P_Bx.
\end{equation}
Applying this in $\bX$ (with $A=\bC$ and $B=\bD$)
and recalling \eqref{e:Pproduct}, we obtain 
\eqref{e:D--Rba}. 
\end{remark}

\begin{remark}[baD--R = CycDyk = CycP in the doubly affine case]
\label{r:doublyaffine}
Let us investigate \eqref{e:wewin3} with starting point $\anfang\in B$ 
further. 
First, we rewrite it as 
\begin{equation}
\label{e:wewin4}
x_0 :=\anfang\in B,\; \text{and}\;
(\forall\knn)\;y_k := P_Bx_k,\;
x_{k+1} := x_k-y_k + P_A(y_k + (\anfang-x_k)/2).
\end{equation}
We now assume additionally that 
\begin{equation}
\text{$A$ and $B$ are linear subspaces of $X$.}
\end{equation}
We claim that for every $k\geq 1$, 
\begin{equation}
\label{e:bond1}
x_k = P_A\anfang - P_BP_A\anfang + P_AP_BP_A\anfang \mp \cdots
+P_A(P_BP_A)^{k-1}\anfang,
\end{equation}
which implies
\begin{equation}
\label{e:bond2}
y_k = (P_BP_A)^k\anfang
\end{equation}
because the differences in \eqref{e:bond1} all lie in $B^\perp$. 
Observe that $y_0=P_Bx_0 = P_B\anfang = \anfang$ since $\anfang\in B$. 
It follows that $x_1 = x_0-y_0 + P_A(y_0+(\anfang-x_0)/2)
= \anfang-\anfang + 
P_A(\anfang + (\anfang-\anfang)/0) = P_A\anfang$ 
which shows that \eqref{e:bond1} holds when $k=1$. 
Now assume that \eqref{e:bond1} holds for some $k\geq 1$.  
Then \eqref{e:bond2} also holds.
Furthermore, 
\begin{equation}
x_k-y_k = P_A\anfang-P_BP_A\anfang\pm \cdots - (P_BP_A)^k\anfang
\end{equation}
and 
\begin{equation}
\anfang - x_k = (\Id-P_A)\anfang + 
(\Id-P_A)P_BP_A\anfang + \cdots (\Id-P_A)(P_BP_A)^{k-1}\anfang
\in A^\perp;
\end{equation}
the latter identity and \eqref{e:bond2} yield
$P_A(y_k+(\anfang-x_k)/2)= P_Ay_k = P_A(P_BP_A)^k\anfang$.
This verifies \eqref{e:bond1} with $k$ replaced by $k+1$;
therefore, by induction, 
\eqref{e:bond1} and hence \eqref{e:bond2} hold true for all $k\geq 1$.
In view of \eqref{e:bond2} and \eqref{e:wewin5}, we observe that
baD--R = CycP = CycDyk in sense that the sequences that arise after
projection onto $B$ are all identical. 
Finally, a translation argument reduces the doubly affine case to the
doubly linear case. 
\end{remark}

\begin{remark}[baD--R $\boldsymbol{\neq}$ CycDyk in general]
\label{r:notequal}
Suppose that $A$ is a nonempty closed convex subset of $X$ and
that $B$ is a subspace of $X$.
Starting both CycDyk and baD--R at $\anfang$, we obtain for $\knn$ 
the update rules 
\begin{equation}
b_0 :=\anfang, p_0 :=0, a_{k+1} := P_A(b_k+p_k), p_{k+1} := b_k+p_k-a_{k+1},
b_{k+1} := P_Ba_{k+1}
\end{equation}
and
\begin{equation}
x_0 := \anfang, y_k := P_Bx_k, x_{k+1} := x_k-y_k+P_A(y_k +
(\anfang-x_k)/2),
\end{equation}
respectively.
One computes that $a_1 = P_A\anfang = x_1$,
$p_1 = \anfang-P_A\anfang$, and 
\begin{equation}
\label{e:wewin6}
b_1 = P_BP_A\anfang = y_1. 
\end{equation}
It thus follows that 
$a_2 = P_A(P_BP_A\anfang + \anfang-P_A\anfang)$
and 
$x_2 = P_{B^\perp}P_A\anfang + P_A(P_BP_A\anfang + 
(\anfang-P_A\anfang)/2))$.
Consequently,
\begin{equation}
\label{e:wewin7}
b_2 = P_BP_A(P_BP_A\anfang + \anfang-P_A\anfang)
\quad\text{and}\quad
y_2 = P_BP_A(P_BP_A\anfang + (\anfang-P_A\anfang)/2)).
\end{equation}
These two vectors certainly appear to be different ---
let us now exhibit a simple example where they actually are different. 
We work in the Euclidean plane and thus assume that $X=\RR^2$.
Set $A := \menge{(\xi,\eta)}{\xi^2 + (\eta-1)^2 \leq 1}$, which 
has the projection formula
\begin{equation}
P_A\colon \RR^2\to\RR^2\colon (\xi,\eta) \mapsto 
\begin{cases}
(\xi,\eta), &\text{if $\xi^2 + (\eta-1)^2 \leq 1$;}\\
\displaystyle (0,1)+ \frac{(\xi,\eta-1)}{\sqrt{\xi^2 + (\eta-1)^2}},
&\text{otherwise,}
\end{cases}
\end{equation}
and set $B := \RR\times \{0\}$. Then
\begin{equation}
P_B\colon\RR^2\to\RR^2\colon (\xi,\eta) \mapsto (\xi,0).
\end{equation}
Set $\anfang=(1,0)$. 
Then \eqref{e:wewin6} turns into 
$b_1 = y_1 = (\sqrt{2}/2,0)$ 
while \eqref{e:wewin7} becomes
\begin{equation}
b_2 = \left( \frac{2}{\sqrt{22-8\sqrt{2}}},0\right)
\approx (0.61180,0) \neq
(0.59718,0) = 
\left( \frac{1}{2}\frac{\sqrt{2}+2}{\sqrt{11-2\sqrt{2}}},0\right)= y_2.
\end{equation}
Let us mention in passing that similar examples exist 
in the product space setting of Section~\ref{ss:mto2} when
some of the constraint sets are balls. 
\end{remark}

\subsection{Summary of best approximation algorithms}

\label{ss:summaryba}

\begin{table*}[H] \centering
\begin{tabular}{@{}lllr@{}} \toprule
Name &Acronym & Formula & Monitor\\ \midrule
Halpern--Wittmann &H--W& \eqref{e:H-W} & $(x_k)_{k\in\NN}$ \\
Cyclic Dykstra &CycDyk& \eqref{e:CycDyk} & $(x_k)_{k\in\NN}$\\
Parallel Dykstra &ParDyk& \eqref{e:ParDyk} & $(x_k)_{k\in\NN}$\\
Haugazeau-like Cyclic Projections &hCycP& \eqref{e:Haug} & $(x_k)_{k\in\NN}$\\
Haugazeau-like Parallel Projections &hParP& \eqref{e:hParP} & $(x_k)_{k\in\NN}$\\
Haugazeau-like Douglas--Rachford&hD--R& \eqref{e:hD--R} & $(\bar{x}_k)_{k\in\NN}$\\
best approximation with Douglas--Rachford&baD--R& \eqref{e:D--Rba} &
$(\bar{x}_k)_{k\in\NN}$\\
\bottomrule
 \end{tabular}
\end{table*}

\section{Nonconvexity}

\label{sec:nonconvexity}

In this section, we provide some remarks on the possible absence of
convexity. 

\subsection{Nonconvex slope constraints}

\label{sec:slope-noncon}
We now consider a case with nonconvex constraints.
This type of constraints does occur frequently in applications;
however, the body of convergence results is sparse and,
to the best of our knowledge, all results are \emph{local}, 
i.e., convergence of the iterates is guaranteed only 
when the starting point is already sufficiently close to the set of 
solutions.
(See, e.g., \cite{BLPW} and the references therein.)
We use the convex constraints from Section~\ref{sec:CPO} with one
crucial modification:
The slope constraints \eqref{e:Si} is tightened by additionally
imposing a \emph{minimum} absolute value slope, i.e., 
we assume the existence of two vectors $a=(\alpha_i)$
and $b=(\beta_i)$ in $\RR^{n-1}_+$ such that 
\begin{equation}
\label{e:Si+}
(\forall i\in I)\quad
S_i := \menge{x=(x_1,\ldots,x_n)\in\RR^n}{
\alpha_i\geq |x_{i+1}-x_i| \geq \beta_i}. 
\end{equation}
Note that if $\beta_i>0$, then $S_i$ is \emph{not convex}. 
Analogously to \eqref{e:Siall},
we aggregate these sets to 
obtain the general nonconvex constraint set
\begin{equation}
S = \bigcap_{i\in\{1,\ldots,n-1\}} S_i.
\end{equation}
Let $i\in\{1,\ldots,n-1\}$ and $x\in X$.
Arguing as in Section~\ref{sss:slopespecial}, 
we obtain the counterpart of \eqref{e:simpleslope} and see that 
$P_{S_i}x$ is different from $x$ at most in positions $i$ and $i+1$,
and\footnote{\ The astute reader will note that when $\beta_i>0$ and 
$x_{i}=x_{i+1}$, then 
$$\big( (P_{S_i}x)_i, (P_{S_i}x)_{i+1}\big) = 
\big\{\thalb(x_i+x_{i+1}+\beta_i,x_i+x_{i+1}-\beta_i),\thalb(x_i+x_{i+1}-\beta_i,x_i+x_{i+1}+\beta_i)
\big\}$$ is not single-valued; indeed, in view of the nonconvexity of $S_i$
and the famous Bunt--Motzkin Theorem 
(see, e.g., \cite[Corollary~21.13]{BC2011}), this is to be expected.
For the actual implementation, one chooses an arbitrary element from 
this set.}
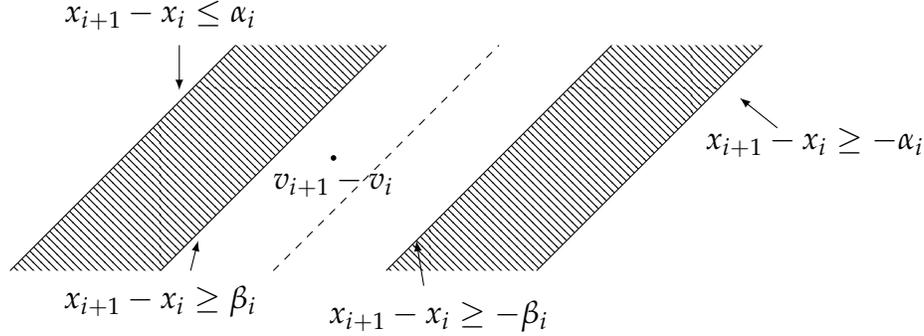
\begin{figure}[H]
\centering
\begin{tikzpicture}
\draw[dashed] (0,0) -- (3,3);
\draw[draw=none,pattern=north west lines] (-3.5,0) -- (-1.5,0) -- (1.5,3) -- (-0.5,3);
\draw[draw=none,pattern=north west lines] (3.5,0) -- (1.5,0) -- (4.5,3) -- (6.5,3);
\draw (-3.5,0) -- (-0.5,3);
\draw (-1.5,0) -- (1.5,3);
\draw (1.5,0) -- (4.5,3);
\draw (3.5,0) -- (6.5,3);
\node[draw=none,fill=none] at (-1.5,3.4) {$x_{i+1}-x_{i} \leq \alpha_{i}$};
\draw[-latex] (-1.25,3) -- (-1.25,2.4);
\node[draw=none,fill=none] at (7.2,1.7) {$x_{i+1}-x_{i} \geq -\alpha_{i}$};
\draw[-latex] (6.7,1.9) -- (6.2,2.3);
\node[draw=none,fill=none] at (-1.5,-0.4) {$x_{i+1}-x_{i} \geq \beta_{i}$};
\draw[-latex] (-1.1,0) -- (-1.0,0.4);
\node[draw=none,fill=none] at (2.2,-0.6) {$x_{i+1}-x_{i} \geq -\beta_{i}$};
\draw[-latex] (2,-0.2) -- (1.9,0.4);
\draw[fill] (0.8,1.5) circle [radius=0.03];
\node [below] at (0.8,1.5) {$\anfang_{i+1}-\anfang_{i}$};
\end{tikzpicture}
\caption{General slope constraint in the nonconvex case.}
\end{figure}

\begin{equation}
\label{e:simpleslope+}
\big( (P_{S_i}x)_i, (P_{S_i}x)_{i+1}\big) = 
\begin{cases}
\thalb(x_i+x_{i+1}+\alpha_i,x_i+x_{i+1}-\alpha_i), &\text{if
$x_{i+1}<x_{i}-\alpha_i$;}\\
(x_i,x_{i+1}), &\text{if $x_i-\alpha_i\leq x_{i+1}\leq x_{i}-\beta_i$;}\\
\thalb(x_i+x_{i+1}+\beta_i,x_i+x_{i+1}-\beta_i), &\text{if
$x_i-\beta_i<x_{i+1}\leq x_i$;}\\
\thalb(x_i+x_{i+1}-\beta_i,x_i+x_{i+1}+\beta_i), &\text{if
$x_i\leq x_{i+1}<x_i+\beta_i$;}\\
(x_i,x_{i+1}), &\text{if $x_i+\beta_i\leq x_{i+1}\leq x_i+\alpha_i$;}\\
\thalb(x_i+x_{i+1}-\alpha_i,x_i+x_{i+1}+\alpha_i), &\text{if
$x_i+\alpha_i<x_{i+1}$.}
\end{cases}
\end{equation}
This formula allows us to deal with the general case as in
Section~\ref{sss:slopegeneral}; thus, we obtain the counterpart
of Proposition~\ref{p:PS}. 

Let us also record the \emph{intrepid counterpart} of
\eqref{e:simpleslope+}:
\begin{equation}
\begin{cases}
\thalb\big(x_i+x_{i+1}+\thalb(\alpha_i+\beta_i),x_i+x_{i+1}-\thalb(\alpha_i+\beta_i)\big),\;
\text{if $x_{i+1}<x_{i}+(\beta_i-3\alpha_i)/2$;}\\[+3mm]
(x_{i+1}+\alpha_i,x_i-\alpha_i),\; \text{if
$x_{i}+(\beta_i-3\alpha_i)/2\leq x_{i+1}<x_i-\alpha_i$;}\\[+3mm]
(x_i,x_{i+1}),\; \text{if $x_i-\alpha_i\leq x_{i+1}\leq
x_{i}-\beta_i$;}\\[+3mm]
(x_{i+1}+\beta_i,x_i-\beta_i),\; \text{if
$x_{i}-\beta_i<x_{i+1}\leq
x_i+\min\{0,(\alpha_i-3\beta_i)/2\}$;}\\[+3mm]
\thalb\big(x_i+x_{i+1}+\thalb(\alpha_i+\beta_i),x_i+x_{i+1}-\thalb(\alpha_i+\beta_i)\big),\;
\text{if $x_i+(\alpha_i-3\beta_i)/2<x_{i+1}\leq x_{i}$;}\\[+3mm]
\thalb\big(x_i+x_{i+1}-\thalb(\alpha_i+\beta_i),x_i+x_{i+1}+\thalb(\alpha_i+\beta_i)\big),\;
\text{if $x_i< x_{i+1}\leq x_i+(3\beta_i-\alpha_i)/2$;}\\[+3mm]
(x_{i+1}-\beta_i,x_i+\beta_i),\; \text{if
$x_i+\max\{0,(3\beta_i-\alpha_i)/2\} \leq x_{i+1} <
x_i+\beta_i$;}\\[+3mm]
(x_i,x_{i+1}),\; \text{if $x_i+\beta_i\leq x_{i+1}\leq
x_{i}+\alpha_i$;}\\[+3mm]
(x_{i+1}-\alpha_i,x_i+\alpha_i),\; \text{if
$x_{i}+\alpha_i < x_{i+1} \leq x_i+(3\alpha_i-\beta_i)/2$;}\\[+3mm]
\thalb\big(x_i+x_{i+1}-\thalb(\alpha_i+\beta_i),x_i+x_{i+1}+\thalb(\alpha_i+\beta_i)\big),\;
\text{if $x_{i}+(3\alpha_i-\beta_i)/2<x_{i+1}$.}
\end{cases}
\end{equation}
There are at least 8 cases; however,
the two cases in the middle arise if and only if $3\beta_i>\alpha_i$.

Having recorded all required formulas, one can now experiment with the
performance of the corresponding algorithms.
In the physics community, Veit Elser has championed
especially D--R with great success 
(see, e.g., \cite{ERT} and \cite{GE}).
Because of its central importance, we consider
in the following subsection a curious cycling behaviour of D--R.

\subsection{Representation of nonconvex constraints and cycling for
D--R}

Suppose momentarily that $C_1,C_2,C_3$ are nonempty closed convex subsets of
$X$ such that $C = C_1\cap C_2\cap C_3\neq\varnothing$. 
Assume further that $C_1\cap C_2$ is still simple enough to admit a
formula for the projector $P_{C_1\cap C_2}$.
The projection methods considered find points in $C$; it does not matter
if we work with either the original three sets $C_1,C_2,C_3$ or
with $C_1\cap C_2$ and $C_3$.
This situation changes dramatically if we are dealing with
\emph{nonconvex} constraints. Performance of projection methods 
crucially depends on the representation of the constraints. 

In fact, the example developed next was discovered numerically. 
It forced us us to deal with nonconvex constraints differently and 
led to the formula \eqref{e:simpleslope+}.

\begin{example}[cycling for D--R]
\label{ex:cycling}
Suppose that $X=\RR^2$, 
that $\alpha\geq\beta>0$, and set
\begin{equation}
C_1 := \menge{(x_1,x_2)}{|x_1-x_2|\leq\alpha}
\quad\text{and}\quad
C_2 := \menge{(x_1,x_2)}{|x_1-x_2|\geq\beta}.
\end{equation}
Clearly, $C_1\cap C_2\neq\varnothing$; in fact, \eqref{e:simpleslope+}
allows us to record a formula for $P_{C_1\cap C_2}$. 
Set 
\begin{equation}
\varepsilon := \beta/5,
\end{equation}
take an arbitrary $\xi\in\RR$, and set
\begin{equation}
\bx_0 = (x_{0,1},x_{0,2}) = 
(\xi,\xi-\varepsilon,\xi-2\varepsilon,\xi+\varepsilon)\in\RR^2\times\RR^2.
\end{equation}
Then the sequences $(\bx_k)_\knn$ and $(\bar{x}_k)_\knn$ generated by
\eqref{e:D--R} do not converge; indeed, they cycle as follows:
\begin{equation}
(\forall\knn)\quad
\bx_{2k} = (\xi,\xi-\varepsilon,\xi-2\varepsilon,\xi+\varepsilon),\;\;
\bx_{2k+1} = (\xi-\varepsilon,\xi,\xi+\varepsilon,\xi-2\varepsilon)
\end{equation}
and
\begin{equation}
(\forall\knn)\quad
\bar{x}_{2k} = (\xi-\varepsilon,\xi),\;\;
\bar{x}_{2k+1} = (\xi,\xi-\varepsilon). 
\end{equation}
\end{example}
\begin{proof}
It is clear that $\bar{x}_0 = (\xi-\varepsilon,\xi)$.
Observe that 
\begin{equation}
\label{e:hawaii1}
(\forall (\xi,\eta)\in\RR^2)\quad 
P_{C_2}(\xi,\eta) = 
\begin{cases}
(\xi,\eta), &\text{if $|\xi-\eta|\geq \beta$;}\\
\thalb(\xi+\eta-\beta,\xi+\eta+\beta), &\text{if $|\xi-\eta|<\beta$ and
$\xi\leq \eta$;}\\
\thalb(\xi+\eta+\beta,\xi+\eta-\beta), &\text{if $|\xi-\eta|<\beta$ and
$\xi\geq \eta$.}
\end{cases}
\end{equation}
By definition,
\begin{subequations}
\begin{align}
x_{1,1} &= x_{0,1} - \bar{x}_0 + P_{C_1}(2\bar{x}_0-x_{0,1})\\
&= (\xi,\xi-\ve)-(\xi-\ve,\xi)+P_{C_1}(2(\xi-\ve,\xi)-(\xi,\xi-\ve))\\
&= (\ve,-\ve)+P_{C_1}(\xi-2\ve,\xi+\ve).
\end{align}
\end{subequations}
Now $|(\xi-2\ve)-(\xi+\ve)| = 3\ve = (3/5)\beta < \beta\leq\alpha$,
so $P_{C_1}$ does not modify $(\xi-2\ve,\xi+\ve)$ and we conclude that 
\begin{equation}
x_{1,1} = (\ve,-\ve)+(\xi-2\ve,\xi+\ve) = (\xi-\ve,\xi).
\end{equation}
Next,
\begin{subequations}
\begin{align}
x_{1,2} &= x_{0,2} - \bar{x}_0 + P_{C_2}(2\bar{x}_0-x_{0,2})\\
&=(\xi-2\ve,\xi+\ve)-(\xi-\ve,\xi)+P_{C_2}(2(\xi-\ve,\xi)-(\xi-2\ve,\xi+\ve))\\
&=(-\ve,\ve)+P_{C_2}(\xi,\xi-\ve).
\end{align}
\end{subequations}
Now $|\xi-(\xi-\ve)| = \ve = \beta/5 < \beta$, so
\eqref{e:hawaii1} yields
\begin{equation}
P_{C_2}(\xi,\xi-\ve) = \thalb(2\xi-\ve+\beta,2\xi-\ve-\beta) 
= (\xi+2\ve,\xi-3\ve).
\end{equation}
It follows that 
$x_{1,2} = (-\ve,\ve) + (\xi+2\ve,\xi-3\ve) = 
(\xi+\ve,\xi-2\ve)$.
Altogether,
\begin{equation}
\bx_1 = (x_{1,1},x_{1,2}) = (\xi-\ve,\xi,\xi+\ve,\xi-2\ve)
\quad\text{and}\quad
\bar{x}_1 = (\xi,\xi-\ve).
\end{equation}
Arguing similarly, we obtain that $\bx_2 = \bx_0$. 
\end{proof}

\begin{remark}
If D--R is started at different points,
as we recommend in Algorithm~\ref{a:D--R}, then
the iterates eventually settle into this cycling behaviour
described in Example~\ref{ex:cycling}. 
\end{remark}

\section{Numerical results}

\label{sec:numeresul}

\subsection{Experimental setup and stopping criteria}

\subsubsection{Generating the set of test problems}
We randomly generate 100 test problems, namely
road splines in groups centered around length 
$L\in\{0.5,1,5,10,20\}$ (unit km).
For each group, we select a design speed
$V\in\{30,50,80,100\}$ (with unit km/h) and maximum elevations
$\xi_{\max}\in\{30,60,100,120,150\}$ (unit m). 
We then generate spline
points $\{(t_1,\anfang_1),\ldots,(t_n,\anfang_n)\}$ such that 
$\anfang\in[0,\xi_{\max}]^n$ and 
$(\forall i\in\{1,\ldots,n-1\})$ 
$\|(t_i,\anfang_i) - (t_{i+1},\anfang_{i+1})\| \geq 0.625 V$,
where $n\in [L/(3\min\{0.625V,30\}),\ldots,1+L/(1.5\min\{0.625V,30\})]\cap\NN$. 
The resulting splines correspond to rather challenging road profiles 
and are therefore ideal for testing. 

\subsubsection{Stopping criteria}
\label{sss:stop}

The constraint sets $C_i$ are as in Section~\ref{sec:CPO} 
(and Section~\ref{sec:slope-noncon} in the nonconvex case). 
Let $(x_k)_{\knn}$ be a sequence to monitor generated by an algorithm
under consideration and set $\varepsilon := 5\cdot 10^{-3}$. 
The feasibility performance measure we employ is 
\begin{equation}
\label{e:d}
d\colon X\to\RP\colon x\mapsto \sqrt{\frac{\sum_{i=1}^{m}
d^2(x,C_i)}{\sum_{i=1}^m d^2(x_0,C_i)}}. 
\end{equation}
Note that $d$ returns the value precisely at points drawn from
$C = C_1\cap\cdots\cap C_m$. 
If we are interested in feasibility only, then we terminate when
\begin{equation}
d(x_k)<\varepsilon;
\end{equation}
in case of best approximation ($P_C\anfang$, where $\anfang\in X$ is given), 
we require additionally that $\|x_k-x_{k-1}\|<\varepsilon$. 
We cap the maximum number of iterations at  $k_{\max} = 5000$.

\subsection{Generation of plots and tables}
Let $\mathcal{P}$ be the set
of test problems and $\mathcal{A}$ be the set of algorithms. 
Let $(x_k^{(a,p)})_\knn$ be the 
monitored sequence generated by algorithm $a\in\mathcal{A}$ applied to the problem
$p\in\mathcal{P}$. 

\subsubsection{Performance plots}
\label{sec:performance}

To compare the performance of the algorithms, we use \emph{performance
profiles}: 
for every $a\in\mathcal{A}$ and 
for every $p\in\mathcal{P}$, we set
\begin{equation}
r_{a,p} := \frac{k_{a,p}}{\min\menge{k_{a',p}}{a'\in\mathcal{A}}} \geq 1,
\end{equation}
where $k_{a,p}\in\{1,2,\ldots,k_{\max}\}$ is the number 
of iterations that $a$ requires to solve $p$ (see
Section~\ref{sss:stop}). 
If $r_{a,p} = 1$, then $a$ uses the least number of iterations to
solve problem $p$. 
If $r_{a,p} > 1$, then $a$ requires $r_{a,p}$
times more iterations for $p$ than the algorithm that uses the least
iterations for $p$.
For each algorithm $a\in\mathcal{A}$, we plot the function
\begin{equation}
\rho_a\colon \RP\to[0,1]\colon \kappa\mapsto 
\frac{\card\menge{p\in\mathcal{P}}{\log_2(r_{a,p})\leq\kappa}}{\card \mathcal{P}},
\end{equation}
where ``$\card$'' denotes the cardinality of a set. 
Thus, $\rho_a(\kappa)$ is the percentage of problems that algorithm
$a$ solves within factor $2^\kappa$ of the best algorithms. 
Therefore, 
an algorithm $a\in\mathcal{A}$ 
is ``fast'' if $\rho_a(\kappa)$ is large for $\kappa$ small;
and $a$ is ``robust'' if $\rho_a(\kappa)$ is large for $\kappa$ large. 
For further information on performance profiles
we refer the reader to \cite{DM}.

\subsubsection{Runtime plots}
\label{sec:approximation}

The sequence $(d(x_k^{(a,p)}))_\knn$ measures the 
runtime progress of $a$ on $p$ with respect to the feasibility performance
measure $d$ (see \eqref{e:d}). 
To get a sense of the average (logarithmatized) progress of each algorithm 
$a\in\mathcal{A}$ at iteration $k$, 
we follow \cite{C97a} and plot the values
of \emph{relative proximity function}, which is defined by 
\begin{subequations}
\begin{align}
\beta_a\colon\NN\to\RR\colon 
k &\mapsto 
10 \log_{10} \Big( \frac{1}{\card\mathcal{P}} 
\sum_{p\in\mathcal{P}} d^2(x_k^{(a,p)})   \Big)\\
&= 
10 \log_{10} \bigg( \frac{1}{\card\mathcal{P}} 
\sum_{p\in\mathcal{P}} \frac{\sum_{i=1}^m \| P_{C_i} x_k^{(a,p)} -
x_k^{(a,p)}\|^2}{\sum_{i=1}^m \| P_{C_i} x_0^{(a,p)} - x_0^{(a,p)}\|^2}\bigg). 
\end{align}
\end{subequations}

\subsubsection{Distance tables}
\label{sec:statsval}

For each algorithm $a\in\mathcal{A}$ and problem $p\in\mathcal{P}$, 
assume that termination occurred at iteration $k(a,p)$. 
We compute the normalized distance to $\anfang$ by 
\begin{equation}
\Delta_{(a,p)} := 
\begin{cases}
\|\anfang - x_{k(a,p)}^{(a,p)}\|\big/\|\anfang\|, &\text{if $k_{(a,p)} < k_{\max}$;}\\
\max \menge{\|\anfang - x_k^{(a',p)}\|}{a'\in\mathcal{A}}\big/\|\anfang\|, &\text{otherwise.}
\end{cases}
\end{equation}
This allows us to consider the collection of normalized distances
$(\Delta_{(a,p)})_{p\in\mathcal{P}}$ for each algorithm $a\in\mathcal{A}$.
Statistical values are recorded for each algorithm in a table allowing us
to compare best approximation performance.

\subsection{Results for feasibility algorithms}

\subsubsection{The convex case}

\label{ss:resfeasconv} 

In this section, we record the results for the convex setting of
Section~\ref{sec:CPO} using the algorithms of Section~\ref{ss:summaryfeas}.

\begin{figure}[H]
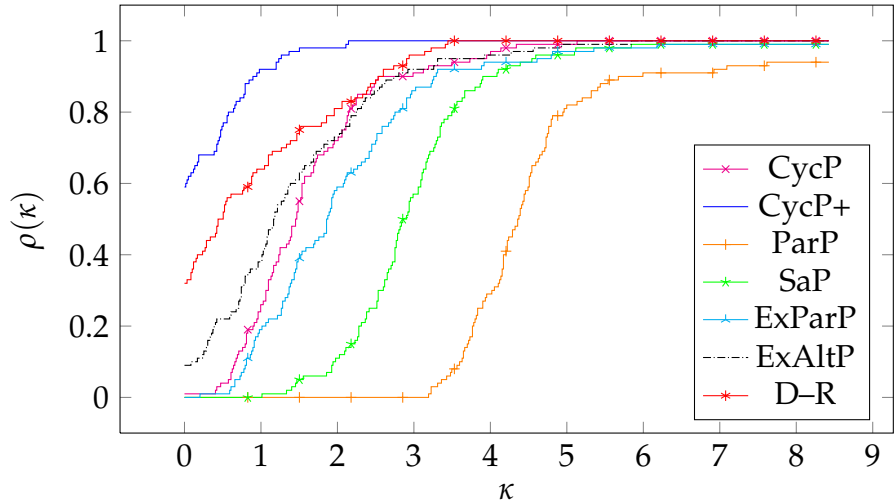

\centering

\caption{Performance profiles for feasibility algorithms in the convex case (see Section~\ref{sec:performance})}
\end{figure}

\begin{figure}[H]
\centering
%
\caption{Runtime plots for feasibility algorithms in the convex case (see Section~\ref{sec:approximation})}
\end{figure}

\begin{table}[H] \centering
\begin{tabular}{@{}lccccccc@{}} \toprule
Algo. & $\Delta$ Min. & $\Delta 1^{\mbox{\tiny st}}$ Qrt. & $\Delta$ Median & $\Delta 3^{\mbox{\tiny rd}}$ Qrt. & $\Delta$ Max. & $\Delta$ Mean & $\Delta$ Std.dev. \\ \midrule
CycP & 0.0563 & 0.3039 & 0.3696 & 0.4806 & 0.9301 & 0.3993 & 0.1568 \\
CycP+ & 0.0563 & 0.2862 & 0.3589 & 0.4403 & 0.7707 & 0.3758 & 0.1399 \\
ParP & 0.0517 & 0.2973 & 0.3687 & 0.4770 & 0.9191 & 0.3990 & 0.1612 \\
SaP & 0.0563 & 0.3033 & 0.3693 & 0.4795 & 0.9292 & 0.3993 & 0.1577 \\
ExParP & 0.0522 & 0.3012 & 0.3697 & 0.4767 & 0.9257 & 0.3981 & 0.1572 \\
ExAltP & 0.0522 & 0.3012 & 0.3697 & 0.4767 & 0.9257 & 0.3973 & 0.1558 \\
D--R & 0.0721 & 0.3833 & 0.4555 & 0.5486 & 1.1353 & 0.4807 & 0.1861 \\
\bottomrule
\end{tabular}
\caption{Statistical data for feasibility algorithms in the convex case 
(see Section~\ref{sec:statsval})}
\label{tab:feastat}
\end{table}

\subsubsection{The nonconvex case}

The plots and tables are analogous to those of
Section~\ref{ss:resfeasconv} except that we work with
the nonconvex slope constraint (see Section~\ref{sec:slope-noncon}).

\begin{figure}[H]
\centering
%
\caption{Performance profiles for feasibility algorithms in the nonconvex case (see Section~\ref{sec:performance})}
\end{figure}

\begin{figure}[H]
\centering
%
\caption{Runtime plots for feasibility algorithms in the nonconvex case (see Section~\ref{sec:approximation})}
\end{figure}

\begin{table}[H] \centering
\begin{tabular}{@{}lccccccc@{}} \toprule
Algo. & $\Delta$ Min. & $\Delta 1^{\mbox{\tiny st}}$ Qrt. & $\Delta$ Median & $\Delta 3^{\mbox{\tiny rd}}$ Qrt. & $\Delta$ Max. & $\Delta$ Mean & $\Delta$ Std.dev. \\ \midrule
CycP & 0.0566 & 0.3039 & 0.3871 & 0.4843 & 0.9301 & 0.4093 & 0.1628 \\
CycP+ & 0.0566 & 0.2520 & 0.3182 & 0.3708 & 0.5873 & 0.3135 & 0.1047 \\
ParP & 0.0522 & 0.3004 & 0.3919 & 0.4853 & 0.9275 & 0.4137 & 0.1653 \\
SaP & 0.0566 & 0.3036 & 0.3820 & 0.4824 & 0.9292 & 0.4027 & 0.1574 \\
ExParP & 0.0527 & 0.3031 & 0.3796 & 0.4786 & 0.9306 & 0.4026 & 0.1603 \\
ExAltP & 0.0526 & 0.3022 & 0.3698 & 0.4771 & 0.9243 & 0.3981 & 0.1566 \\
D--R & 0.0756 & 0.3835 & 0.4567 & 0.5486 & 1.1044 & 0.4810 & 0.1846 \\
\bottomrule
\end{tabular}
\caption{Statistical data for feasibility algorithms in the nonconvex case 
(see Section~\ref{sec:statsval})}
\label{tab:nc-feasstat}
\end{table}

\subsubsection{Conclusions}

Overall,  CycP+, ExAltP, and D--R emerge as good algorithms for
feasibility. 
CycP+ yields solutions closest to the initial spline and is among the most
robust. 
In the convex case, CycP+ is also the fastest algorithm;
in the nonconvex case, D--R is faster. 
Both ExAltP and D--R are good choices when the number of iterations 
is small.

\subsection{Results for superiorized feasibility algorithms}

\label{ss:ressuper}

\subsubsection{The convex case}

\label{sss:ressuperconv}

In this section, we record the results for the convex setting of
Section~\ref{sec:CPO} using superiorized feasibility algorithms (see
Section~\ref{ss:superior}). 

\begin{figure}[H]
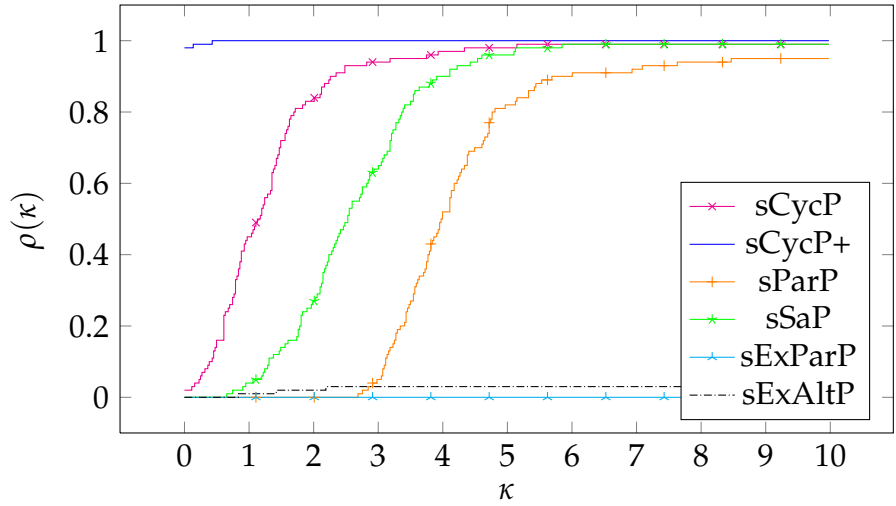

\centering

\caption{Performance profiles for superiorized feasibility algorithms in
the convex case (see Section~\ref{sec:performance})}
\end{figure}

\begin{figure}[H]
\centering
%
\caption{Runtime plots for superiorized feasibility algorithms in the
convex case (see Section~\ref{sec:approximation})}
\end{figure}

\begin{table}[H] \centering
\begin{tabular}{@{}lccccccc@{}} \toprule
Algo. & $\Delta$ Min. & $\Delta 1^{\mbox{\tiny st}}$ Qrt. & $\Delta$ Median & $\Delta 3^{\mbox{\tiny rd}}$ Qrt. & $\Delta$ Max. & $\Delta$ Mean & $\Delta$ Std.dev. \\ \midrule
sCycP & 0.0519 & 0.3039 & 0.3696 & 0.4805 & 0.9301 & 0.3992 & 0.1570 \\
sCycP+ & 0.0563 & 0.2863 & 0.3599 & 0.4410 & 0.7707 & 0.3761 & 0.1399 \\
sParP & 0.0517 & 0.2973 & 0.3659 & 0.4762 & 0.9191 & 0.3940 & 0.1559 \\
sExParP & 0.0563 & 0.3039 & 0.3711 & 0.4825 & 0.9301 & 0.4018 & 0.1559 \\
sExAltP & 0.0519 & 0.3039 & 0.3711 & 0.4825 & 0.9301 & 0.4017 & 0.1561 \\
sSaP & 0.0519 & 0.3033 & 0.3693 & 0.4795 & 0.9292 & 0.3985 & 0.1567 \\
\bottomrule
\end{tabular}\caption{Statistical data for superiorized feasibility algorithms in the
convex case (see Section~\ref{sec:statsval})}
\label{tab:supdev}
\end{table}

\subsubsection{The nonconvex case}

\label{sss:ressupernonconv}

The plots and tables are analogous to those of
Section~\ref{sss:ressuperconv} except that we work with
the nonconvex slope constraint (see Section~\ref{sec:slope-noncon}).

\begin{figure}[H]
\centering
%
\caption{Performance profiles for superiorized feasibility algorithms in
the nonconvex case (see Section~\ref{sec:performance})}
\end{figure}

\begin{figure}[H]
\centering
%
\caption{Runtime plots for superiorized feasibility algorithms in the
nonconvex case (see Section~\ref{sec:approximation})}
\end{figure}

\begin{table}[H] \centering
\begin{tabular}{@{}lccccccc@{}} \toprule
Algo. & $\Delta$ Min. & $\Delta 1^{\mbox{\tiny st}}$ Qrt. & $\Delta$ Median & $\Delta 3^{\mbox{\tiny rd}}$ Qrt. & $\Delta$ Max. & $\Delta$ Mean & $\Delta$ Std.dev. \\ \midrule
sCycP & 0.0566 & 0.3039 & 0.3754 & 0.4759 & 0.9301 & 0.4002 & 0.1563 \\
sCycP+ & 0.0566 & 0.2885 & 0.3594 & 0.4388 & 0.7708 & 0.3765 & 0.1392 \\
sParP & 0.0522 & 0.3004 & 0.3754 & 0.4744 & 0.9275 & 0.3973 & 0.1554 \\
sExAltP & 0.0523 & 0.3042 & 0.3754 & 0.4762 & 0.9301 & 0.4030 & 0.1555 \\
sExParP & 0.0566 & 0.3042 & 0.3754 & 0.4762 & 0.9301 & 0.4030 & 0.1554 \\
sSap & 0.0522 & 0.3036 & 0.3714 & 0.4734 & 0.9292 & 0.3978 & 0.1556 \\
\bottomrule
\end{tabular}\label{tab:nc-supstat}
\caption{Statistical data for superiorized feasibility algorithms in the
nonconvex case (see Section~\ref{sec:statsval})}
\end{table}

\subsubsection{Conclusions}

We clearly see that sCycP+ is not only the fastest but also the most 
robust superiorized feasibility algorithm.

\subsection{Results for best approximation algorithms}

\subsubsection{The convex case}

\label{sss:resbaconv}

In this section, we record the results for the convex setting of
Section~\ref{sec:CPO} using the best approximation algorithms of
Section~\ref{sec:ba}.

\begin{figure}[H]
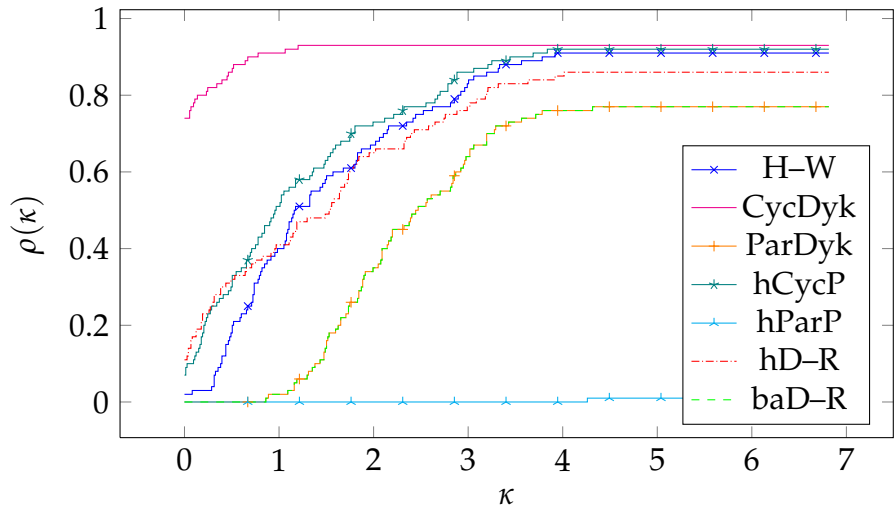

\centering

\caption{Performance profiles for best approximation algorithms in the
convex case (see Section~\ref{sec:performance})}
\end{figure}

\begin{figure}[H]
\centering
%
\caption{Runtime plots for best approximation algorithms in the convex case
(see Section~\ref{sec:approximation})}
\end{figure}

\begin{table}[H] \centering
\begin{tabular}{@{}lccccccc@{}} \toprule
Algo. & $\Delta$ Min. & $\Delta 1^{\mbox{\tiny st}}$ Qrt. & $\Delta$ Median & $\Delta 3^{\mbox{\tiny rd}}$ Qrt. & $\Delta$ Max. & $\Delta$ Mean & $\Delta$ Std.dev. \\ \midrule
H--W & 0.0525 & 0.2978 & 0.3644 & 0.4800 & 0.9251 & 0.3959 & 0.1570 \\
CycDyk & 0.0519 & 0.3022 & 0.3695 & 0.4862 & 0.9379 & 0.4006 & 0.1584 \\
ParDyk & 0.0517 & 0.2962 & 0.3656 & 0.4799 & 0.9379 & 0.3972 & 0.1576 \\
hCycP & 0.0519 & 0.2986 & 0.3645 & 0.4801 & 0.9251 & 0.3962 & 0.1569 \\
hParP & 0.0517 & 0.3042 & 0.3720 & 0.4864 & 0.9379 & 0.4016 & 0.1585 \\
hD--R & 0.0519 & 0.3042 & 0.3720 & 0.4864 & 0.9337 & 0.4013 & 0.1579 \\
D--Rba & 0.0517 & 0.2962 & 0.3656 & 0.4799 & 0.9379 & 0.3972 & 0.1576 \\
\bottomrule
\end{tabular}
\caption{Statistical data for best approximation algorithms in the convex
case (see Section~\ref{sec:statsval})}
\label{tab:supstat}
\end{table}

\subsubsection{The nonconvex case}

\label{sss:resbanonconv}

The plots and tables are analogous to those of
Section~\ref{sss:resbaconv} except that we work with
the nonconvex slope constraint (see Section~\ref{sec:slope-noncon}).

\begin{figure}[H]
\centering
%
\caption{Performance profiles for best approximation algorithms in the
nonconvex case (see Section~\ref{sec:performance})}
\end{figure}

\begin{figure}[H]
\centering
%
\caption{Runtime plots for best approximation algorithms in the nonconvex case
(see Section~\ref{sec:approximation})}
\end{figure}

\begin{table}[H] \centering
\begin{tabular}{@{}lccccccc@{}} \toprule
Algo. & $\Delta$ Min. & $\Delta 1^{\mbox{\tiny st}}$ Qrt. & $\Delta$ Median & $\Delta 3^{\mbox{\tiny rd}}$ Qrt. & $\Delta$ Max. & $\Delta$ Mean & $\Delta$ Std.dev. \\ \midrule
H--W & 0.0525 & 0.3019 & 0.3794 & 0.4825 & 0.9252 & 0.3998 & 0.1570 \\
CycDyk & 0.0667 & 0.3053 & 0.3861 & 0.4863 & 0.9385 & 0.4079 & 0.1539 \\
ParDyk & 0.0667 & 0.3006 & 0.3846 & 0.4863 & 0.9385 & 0.4060 & 0.1529 \\
hCycP & 0.0522 & 0.3022 & 0.3782 & 0.4808 & 0.9253 & 0.3998 & 0.1568 \\
hParP & 0.0520 & 0.3057 & 0.3863 & 0.4864 & 0.9385 & 0.4069 & 0.1570 \\
hD--R & 0.0652 & 0.3057 & 0.3825 & 0.4864 & 0.9342 & 0.4074 & 0.1540 \\
baD--R & 0.0667 & 0.3006 & 0.3846 & 0.4863 & 0.9385 & 0.4060 & 0.1529 \\
\bottomrule
\end{tabular}
\caption{Statistical data for best approximation algorithms in the nonconvex
case (see Section~\ref{sec:statsval})}
\label{tab:nc-bapstat}
\end{table}

\subsubsection{ParDyk vs baD--R}

The plots and tables of Sections~\ref{sss:resbaconv}
and \ref{sss:resbanonconv} 
suggest that ParDyk and baD--R are the same algorithms. 
This is not true for these algorithms in their full generality
as we pointed out in Remark~\ref{r:notequal}. 
A possible explanation for this identical performance could be 
that the algorithms produce sequences that stay either outside or
inside the halfspaces comprising the constraints. 
If that is indeed the case, then
the projection onto the halfspace would be indistinguishable from
the projection onto either a hyperplane or the entire space;
consequently, the iterates would be identical by
Remark~\ref{r:doublyaffine}.

\subsubsection{Conclusions}

hCycP is a good robust choice for both convex and nonconvex problems;
CycDyk does well for convex problems, and H--W does well for nonconvex
problems. 

Examining the data for all three types of algorithms
(feasibility, superiorization, best approximation), we can make
a compelling case for CycP+, which emerges as the best overall algorithm with
respect to speed, robustness, and best approximation properties.

\subsection{Projection methods vs Linear Programming (LP) algorithms}

\label{sec:lpsimplex}

In Section \ref{sec:swissarmy}, we made the case for projection methods and
commented on their competitiveness with other optimization methods. 
To illustrate this claim, we interpreted \eqref{e:www}
as the constraints of a \emph{Linear Programming (LP)} problem.
We then solved our test problems with the
\emph{GNU Linear Programming Kit (GLPK)} \cite{GLPK}.
Two objective functions were employed:
$x\mapsto 0$, which corresponds to the feasibility problem \eqref{e:cfp},
and $x\mapsto \|x-\anfang\|_1$, which is similar to
the best approximation problem\footnote{ Here, $\|x\|_1 =
\sum_{i=1}^{n}|x_i|$ denotes the $\ell_1$-norm of $x\in X$.}
\eqref{e:130111a}. We call these methods GLPK0 and GLPK1, respectively \cite[page~257]{Williams}. 
We shall compare these two methods against the overall best projection
method, CycP+.

\begin{table*}[H] \centering

\caption{Performance profiles for convex feasibility problems 
comparing the number of iterations of CycP+, GLPK0 and GLPK1 until
convergence (see Section~\ref{sec:performance})}
\label{fig:lp}
\end{figure}

\begin{table}[H] \centering
\begin{tabular}{@{}lccccccc@{}} \toprule
Algo. & $\Delta$ Min. & $\Delta 1^{\mbox{\tiny st}}$ Qrt. & $\Delta$ Median & $\Delta 3^{\mbox{\tiny rd}}$ Qrt. & $\Delta$ Max. & $\Delta$ Mean & $\Delta$ Std.dev. \\ \midrule
CycP+ & 0.0563 & 0.2862 & 0.3589 & 0.4403 & 0.7707 & 0.3758 & 0.1399 \\
GLPK0 & 0.1938 & 0.6636 & 0.8030 & 0.9281 & 1.0039 & 0.7767 & 0.1813 \\
GLPK1 & 0.0000 & 0.3062 & 0.3813 & 0.5001 & 0.9406 & 0.4076 & 0.1625 \\
\bottomrule
\end{tabular}
\caption{Statistical data in the convex case (see
Section~\ref{sec:statsval})}
\label{tab:avpstat}
\end{table}

The nonconvex slope constraints from Section~\ref{sec:slope-noncon}
can be handled by using \emph{binary variables}; this leads
to a \emph{Mixed Integer Linear Programming (MIP)} problem (see, e.g.,
\cite[Chapter~12]{HL}), which GLPK \cite{GLPK} can handle as well.

\begin{figure}[H]
\centering
%
\caption{Performance profiles for nonconvex feasibility problems comparing
the number of iterations of CycP+, GLPK0 and GLPK1 until convergence (see
Section~\ref{sec:performance})}
\label{fig:mip}
\end{figure}

\begin{table}[H] \centering
\begin{tabular}{@{}lccccccc@{}} \toprule
Algo. & $\Delta$ Min. & $\Delta 1^{\mbox{\tiny st}}$ Qrt. & $\Delta$ Median & $\Delta 3^{\mbox{\tiny rd}}$ Qrt. & $\Delta$ Max. & $\Delta$ Mean & $\Delta$ Std.dev. \\ \midrule
CycP+ & 0.0566 & 0.2520 & 0.3182 & 0.3708 & 0.5873 & 0.3135 & 0.1047 \\
GLPK0 & 0.1938 & 0.5371 & 0.6245 & 0.8240 & 6.2070 & 0.8226 & 0.7464 \\
GLPK1 & 0.0840 & 0.3127 & 0.3847 & 0.5001 & 2.1468 & 0.4261 & 0.2325 \\
\bottomrule
\end{tabular}
\caption{Statistical data in the nonconvex case (see
Section~\ref{sec:statsval})}
\label{tab:nc-avpstat}
\end{table}

\section{Concluding remarks}
\label{sec:conclusion}

Using the practical example of road design, we formulated 
(convex and nonconvex) feasibility and best approximation problems. 
We studied projection methods and implemented them to solve the
feasibility problems. A clear winner emerged, even when compared to LP
solvers: CycP+, the method of cyclic intrepid projections. 

In the future, we plan to study the influence of parameters on the
performance of the algorithms presented. The design and testing of hybrid methods
that aim to combine advantageous traits of various algorithms is also of
considerable interest. 
Finally, rigorous convergence statements about CycP+ in the nonconvex setting await 
to be discovered.

\section*{Acknowledgment}
HHB was partially supported by the Natural Sciences and
Engineering Research Council of Canada and by the Canada Research Chair
Program.


\begin{thebibliography}{999}


\bibitem{AASHTO}
American Association of State Highway and Transportation
Officials, 
\emph{A Policy on Geometric Design of Highways and Streets},
sixth edition, 
Washington, D.C., 2011.  


\bibitem{BBDyk}
H.H.\ Bauschke and J.M.\ Borwein,
Dykstra's alternating projection algorithm for two sets,
\emph{Journal of Approximation Theory}~79 (1994), 418--443. 

\bibitem{BB96}
H.H.\ Bauschke and J.M.\ Borwein,
On projection algorithms for solving convex feasibility problems,
\emph{SIAM Review}~38(3) (1996), 367--426. 

\bibitem{BBHM}
H.H.\ Bauschke, R.I.\ Bo\c{t}, W.L.\ Hare, and W.M.\ Moursi,
Attouch-Th\'era duality revisited: paramonotonicity and operator
splitting,
\emph{Journal of Approximation Theory}~164 (2012), 1065--1084. 

\bibitem{MOR}
H.H.\ Bauschke and P.L.\ Combettes,
A weak-to-strong convergence principle for
Fej\'er-monotone methods in Hilbert spaces,
\emph{Mathematics of Operations Research}~26(2) (2001), 248--264. 

\bibitem{BC2011}
H.H.\ Bauschke and P.L.\ Combettes,
\emph{Convex Analysis and Monotone Operator Theory in Hilbert Spaces},
Springer, 2011.

\bibitem{BCK}
H.H.\ Bauschke, P.L.\ Combettes, and S.G.\ Kruk, 
Extrapolation algorithm for affine-convex feasibility problems, 
\emph{Numerical Algorithms}~41(3) (2006), 239--274.

\bibitem{BCL}
H.H.\ Bauschke, P.L.\ Combettes, and D.R.\ Luke,
Finding best approximation pairs relative to two
closed convex sets in Hilbert spaces,
\emph{Journal of Approximation Theory}~127 (2004), 178--192. 

\bibitem{BK}
H.H.\ Bauschke and S.G.\ Kruk,
Reflection-projection method for convex feasibility problems
with an obtuse cone,
\emph{Journal of Optimization Theory and Applications}~12(3) (2004),
503--531. 

\bibitem{BLPW}
H.H.\ Bauschke, D.R.\ Luke, H.M.\ Phan, and X.\ Wang,
Restricted normal cones and the method of alternating projections,
\texttt{http://arxiv.org/abs/1205.0318v1}, May 2012.

\bibitem{Bellman}
R.\ Bellman,
On the approximation of curves by line segments using dynamic programming,
\emph{Communications of the ACM}~4(6) (1961), 284. 

\bibitem{Boyd}
S.\ Boyd, N.\ Parikh, E.\ Chu, B.\ Peleato, and J.\ Eckstein,
Distributed optimization and statistical learning via
the alternating direction method of multipliers,
\emph{Foundations and Trends in Machine Learning}~3(1) (2010), 1--122. 

\bibitem{Boyd2012}
S.\ Boyd and J.\ Mattingley,
CVXGEN: A Code Generator for Embedded Convex Optimization,
\emph{Optimization and Engineering}~13(1) (2012), 1--27. 

\bibitem{B-D}
J.P.\ Boyle and R.L.\ Dykstra,
A method for finding projections onto the intersection of convex sets
in Hilbert spaces, in
\emph{Advances in Order Restricted Statistical Interference}, Lecture
Notes in Statistics~37, Springer, 1986. 

\bibitem{Cegielski}
A.\ Cegielski,
\emph{Iterative Methods for Fixed Point Problems in Hilbert Spaces},
Springer, 2012. 

\bibitem{CCCDH}
Y.\ Censor, W.\ Chen, P.L.\ Combettes, R.\ Davidi, and G.T.\ Herman,
On the effectiveness of projection methods for convex feasibility
problems with linear inequality constraints, 
\emph{Computational Optimization and Applications}~51(3) (2012),
1065--1088.

\bibitem{CCP}
Y.\ Censor, W.\ Chen, and H.\ Pajoohesh,
Finite convergence of a subgradient projections method
with expanding controls,
\emph{Applied Mathematics and Optimization}~64 (2011), 273--285. 

\bibitem{CDH}
Y.\ Censor, R.\ Davidi, and G.T.\ Herman,
Pertubation resilience and superiorization of iterative algorithms,
\emph{Inverse Problems}~26(6) (2010), 065008 (12 pages). 

\bibitem{CE}
Y.\ Censor and T.\ Elfving,
{New methods for linear inequalities},
\emph{Linear Algebra and its Applications}~42 (1982), 199--211. 

\bibitem{CEH}
Y.\ Censor, T.\ Elfving, and G.T.\ Herman,
Averaging strings of sequential iterations for convex
feasibility problems, in 
\emph{Inherently Parallel Algorithms in Feasibility and 
Optimization and Their Applications},
D.~Butnariu, Y.~Censor, and S.~Reich (editors),
pp.~101--114, Elsevier, 2001. 

\bibitem{CT}
Y.\ Censor and E.\ Tom,
Convergence of string averaging projection schemes for 
inconsistent feasibility problems,
\emph{Optimization Methods and Software}~18 (2003), 543--554.

\bibitem{CZ}
Y.\ Censor and S.A\ Zenios,
\emph{Parallel Optimization},
Oxford University Press, 1997. 

\bibitem{C97a}
P.L.\ Combettes,
Convex set theoretic image recovery
by extrapolated iterations of parallel subgradient
projections, 
\emph{IEEE Transactions on Image Processing}~6(4) (1997), 493--506. 

\bibitem{C97}
P.L.\ Combettes,
Hilbertian convex feasibility problems: 
convergence of projection methods,
\emph{Applied Mathematics \& Optimization}~35 (1997), 311--330. 

\bibitem{CDV}
P.L.\ Combettes, {D{\hspace{-1.6ex}{\raise 0.4ex\hbox{-}\hspace{.8ex}}}inh~D\~ung}, and B.C.\ V\~u, 
Dualization of signal recovery problems, 
\emph{Set--Valued and Variational Analysis}~18 (2010), 373--404.

\bibitem{CP}
P.L.\ Combettes and J.-C.\ Pesquet,
Proximal splitting methods in signal processing,
Chapter~10 in \emph{Fixed-Point Algorithms for Inverse Problems in
Science and Engineering}, Springer, 2011. 

\bibitem{deBoor}
C.\ de Boor,
\emph{A Practical Guide to Splines},
revised edition, 
Springer, 2001. 

\bibitem{Deutsch}
F.\ Deutsch,
\emph{Best Approximation in Inner Product Spaces}, 
Springer, 2001. 

\bibitem{DM}
E.D.\ Dolan and J.J.\ Mor\'{e},
Benchmarking optimization software with performance profiles,
\emph{Mathematical Programming (Series A)} 91 (2002), 201--213.

\bibitem{EckBer}
J.\ Eckstein and D.P.\ Bertsekas,
\emph{On the Douglas-Rachford splitting method
and the proximal point algorithm for maximal monotone
operators},
\emph{Mathematical Programming (Series A)} 55 (1992), 293--318.

\bibitem{ERT}
V.\ Elser, I.\ Rankenburg, and P.\ Thibault,
Searching with iterated maps,
\emph{Proceedings of the National Academy of Sciences}~104(2) (2007),
418--423. 

\bibitem{GM}
D.\ Gabay and B.\ Mercier,
A dual algorithm for the solution of nonlinear variational problems
via finite elements approximations,
\emph{Compututers \& Mathematics with Applications}~2 (1976), 17--40.

\bibitem{GlM}
R.\ Glowinski and A.\ Marrocco,
Sur l'approximation, par \'el\'ements finis d'ordre un, et la
r\'esolution, par p\'enalisation-dualit\'e, d'une classe 
de probl\`emes de Dirichlet non lin\'eaires,
\emph{RAIRO Analyse Num\'erique}~2 (1975), 41--76. 


\bibitem{GLPK}
GNU Linear Programming Kit, Version 4.45, 2010. Available at URL 
\url{http://www.gnu.org/software/glpk/glpk.html}


\bibitem{GK}
K.\ Goebel and W.A.\ Kirk,
\emph{Topics in Metric Fixed Point Theory},
Cambridge University Press, 1990.

\bibitem{GR}
K.\ Goebel and S.\ Reich,
\emph{Uniform Convexity, Hyperbolic Geometry, and Nonexpansive Mappings},
Marcel Dekker, 1984.

\bibitem{Gould}
N.I.M.\ Gould,
How good are projection methods for convex feasibility problems?
\emph{Computational Optimization and Applications}~40 (2008), 1--12. 

\bibitem{GE}
S.\ Gravel and V.\ Elser,
Divide and concur: a general approach to constraint satisfaction,
\emph{Physical Review E}~78 (2008) 036706 (5 pages). 

\bibitem{Halpern}
B.\ Halpern,
Fixed points of nonexpanding maps,
\emph{Bulletin of the AMS}~73 (1967), 957--961.

\bibitem{Haug}
Y.\ Haugazeau,
\emph{Sur les In\'equations Variationnelles et la
Minimisation de Foctionnelles Convexes},
Ph.D.\ thesis, Universit\'e de Paris, 1968. 

\bibitem{Herman}
G.T.\ Herman,
\emph{Fundamentals of Computerized Tomography},
second edition, Springer, 2009.

\bibitem{HC}
G.T.\ Herman and W.\ Chen,
A fast algorithm for solving a linear feasibility problem
with application to intensity-modulated radiation therapy,
\emph{Linear Algebra and its Applications}~428 (2008), 1207--1217. 

\bibitem{HL}
F.S.\ Hillier and G.J.\ Lieberman,
\emph{Introduction to Operations Research},
seventh edition, 
McGraw Hill, 2001. 

\bibitem{LM}
P.-L.\ Lions and B.\ Mercier,
Splitting algorithms for the sum of two nonlinear
operators,
\emph{SIAM Journal on Numerical Analysis}~16 (1979), 964--979.

\bibitem{Koch}
V.R.\ Koch,
Road Design Optimization,
\emph{US patent application no.} 13/626,451, filed on September~25, 2012.

\bibitem{Reich83}
S.\ Reich,
A limit theorem for projections,
\emph{Linear and Multilinear Algebra}~13 (1983), 281--290.

\bibitem{Schumaker}
L.L.\ Schumaker,
\emph{Spline Functions: Basic Theory}, 
third edition, 
Cambridge University Press, 2007.

\bibitem{Williams}
H.P.\ Williams,
\emph{Model Building in Mathematical Programming},
third edition,
Wiley,1990.

\bibitem{Wittmann}
R.\ Wittmann,
Approximation of fixed points of nonexpansive mappings,
\emph{Archiv der Mathematik}~58 (1992), 486--491.


\end{thebibliography}
\end{document}